\documentclass{amsart}
\begin{document}
\title[Eigenforms]{Hecke operators on differential modular forms mod $p$}
\author{Alexandru Buium and Arnab Saha}
\def \dpi{\d_{\pi}}
\def \bT{{\bf T}}
\def \cI{{\mathcal I}}
\def \cJ{{\mathcal J}}
\def \ZN{\bZ[1/N,\zeta_N]}
\def \tA{\tilde{A}}
\def \o{\omega}
\def \tB{\tilde{B}}
\def \tC{\tilde{C}}
\def \alph{A}
\def \bet{B}
\def \bsigma{\bar{\sigma}}
\def \y{^{\infty}}
\def \Ra{\Rightarrow}
\def \uBS{\overline{BS}}
\def \lBS{\underline{BS}}
\def \lB{\underline{B}}
\def \<{\langle}
\def \>{\rangle}
\def \hL{\hat{L}}
\def \cU{\mathcal U}
\def \cF{\mathcal F}
\def \S{\Sigma}
\def \st{\stackrel}
\def \sd{Spec_{\d}\ }
\def \pd{Proj_{\d}\ }
\def \s{\sigma_2}
\def \i{\sigma_1}
\def \bs{\bigskip}
\def \cD{\mathcal D}
\def \cC{\mathcal C}
\def \cT{\mathcal T}
\def \cK{\mathcal K}
\def \cX{\mathcal X}
\def \sX{X_{set}}
\def \cY{\mathcal Y}
\def \cS{X}
\def \cR{\mathcal R}
\def \cE{\mathcal E}
\def \tcE{\tilde{\mathcal E}}
\def \cP{\mathcal P}
\def \cA{\mathcal A}
\def \cV{\mathcal V}
\def \cM{\mathcal M}
\def \cL{\mathcal L}
\def \cN{\mathcal N}
\def \tcM{\tilde{\mathcal M}}
\def \caS{\mathcal S}
\def \cG{\mathcal G}
\def \cB{\mathcal B}
\def \tG{\tilde{G}}
\def \cF{\mathcal F}
\def \h{\hat{\ }}
\def \hp{\hat{\ }}
\def \tS{\tilde{S}}
\def \tP{\tilde{P}}
\def \tA{\tilde{A}}
\def \tX{\tilde{X}}
\def \tcS{\tilde{X}}
\def \tT{\tilde{T}}
\def \tE{\tilde{E}}
\def \tV{\tilde{V}}
\def \tC{\tilde{C}}
\def \tI{\tilde{I}}
\def \tU{\tilde{U}}
\def \tG{\tilde{G}}
\def \tu{\tilde{u}}
\def \chu{\check{u}}
\def \tx{\tilde{x}}
\def \tL{\tilde{L}}
\def \tY{\tilde{Y}}
\def \d{\delta}
\def \e{\chi}
\def \bW{\mathbb W}
\def \bV{{\mathbb V}}
\def \bF{{\bf F}}
\def \bE{{\bf E}}
\def \bC{{\bf C}}
\def \bO{{\bf O}}
\def \bR{{\bf R}}
\def \bA{{\bf A}}
\def \bB{{\bf B}}
\def \cO{\mathcal O}
\def \ra{\rightarrow}
\def \bx{{\bf x}}
\def \f{{\bf f}}
\def \bX{{\bf X}}
\def \bH{{\bf H}}
\def \bS{{\bf S}}
\def \bF{{\bf F}}
\def \bN{{\bf N}}
\def \bK{{\bf K}}
\def \bE{{\bf E}}
\def \bB{{\bf B}}
\def \bQ{{\bf Q}}
\def \bd{{\bf d}}
\def \bY{{\bf Y}}
\def \bU{{\bf U}}
\def \bL{{\bf L}}
\def \bQ{{\bf Q}}
\def \bP{{\bf P}}
\def \bR{{\bf R}}
\def \bC{{\bf C}}
\def \bD{{\bf D}}
\def \bM{{\bf M}}
\def \bZ{{\mathbb Z}}
\def \xtoleqr{x^{(\leq r)}}
\def \hU{\hat{U}}
\def \k{\kappa}
\def \ee{\overline{p^{\k}}}

\newtheorem{THM}{{\!}}[section]
\newtheorem{THMX}{{\!}}
\renewcommand{\theTHMX}{}
\newtheorem{theorem}{Theorem}[section]
\newtheorem{corollary}[theorem]{Corollary}
\newtheorem{lemma}[theorem]{Lemma}
\newtheorem{proposition}[theorem]{Proposition}
\theoremstyle{definition}
\newtheorem{definition}[theorem]{Definition}
\theoremstyle{remark}
\newtheorem{remark}[theorem]{Remark}
\newtheorem{example}[theorem]{\bf Example}
\numberwithin{equation}{section}
\address{University of New Mexico \\ Albuquerque, NM 87131}
\email{buium@math.unm.edu, arnab@math.unm.edu} \subjclass[2000]{11 F 32, 11 F 85}
\maketitle

\begin{abstract}
A   description is given
of all primitive   $\d$-series mod $p$ of order $1$  which are eigenvectors of all the Hecke operators $nT_{\k}(n)$, ``$pT_{\k}(p)$", $(n,p)=1$, and which are $\d$-Fourier expansions of $\d$-modular forms of arbitrary order and weight $w$ with $deg(w)=\k\geq 0$; this set of $\d$-series is shown to be in a natural one-to-one correspondence with the set of series mod $p$ (of order $0$) which are eigenvectors of all the Hecke operators $T_{\k+2}(n)$, $T_{\k+2}(p)$, $(n,p)=1$ and which are Fourier expansions of (classical) modular forms of weight $\equiv \k+2$ mod $p-1$.
\end{abstract}

\section{Introduction}

This present  paper is a direct continuation of \cite{eigen} and , implicitly, of a series of papers devoted to the study of arithmetic differential equations  \cite{char, difmod, Barcau, book,  dcc,igusa}; however, for the convenience of the reader,
the present paper is written so as to be logically independent of \cite{eigen} and of the other above cited papers. Rather, we will quickly review the relevant material from some of these papers  as needed.

The plan of  this Introduction is as follows. We begin by  quickly recalling the basic definitions of this theory following \cite{char,book}.
For more details on some of these definitions we refer to the body of the present paper. Then we will state our main result (Theorem \ref{maint}). Finally we will make some comments on the larger picture and motivations beyond this theory.

\subsection{$\d$-functions \cite{char,book}}
A map $\d:A \ra B$ from a ring $A$ into a  $p$-torsion free $A$-algebra $B$
 is called a {\it $p$-derivation} if the map $\phi:A \ra B$, $\phi(x)=x^p+p\d x$, is a ring homomorphism. When $\d$ is given $\phi$ will always have the meaning above. A ring equipped with a $p$-derivation will be refered to as a {\it $\d$-ring}.
 Denote by $R$ the completion of the maximum unramified extension of the ring of $p$-adic integers. Set $k=R/pR$, $K=R[1/p]$, let $\phi:R \ra R$ be the unique ring automorphism lifting the $p$-power Frobenius $F:k \ra k$, and denote by $\d:R \ra R$ the $p$-derivation
$\d x=\frac{\phi(x)-x^p}{p}$. This makes $R$ a $\d$-ring and this $\d$-ring structure on $R$ is unique.
For any affine smooth scheme $V \subset {\mathbb A}^m$ over $R$ a function $f:V(R)\ra R$ will be called a {\it $\d$-function
of order $r$} on $V$ \cite{char} if
there exists a restricted power series $\Phi$ in $m(r+1)$ variables, with $R$-coefficients such that
$f(a)=\Phi(a,\d a,...,\d^r a)$,
for all $a\in V(R)\subset R^m$.
(Recall that a power series is called {\it restricted} is its coefficients tend to $0$ $p$-adically.)
We denote by $\cO^r(V)$ the ring of $\d$-functions of order $r$ on $V$.
This concept can be naturally extended to the non-affine case \cite{char} but we will not need this
extension in the present paper.

\subsection{$\d$-modular forms \cite{difmod,book}}
Let $N>4$ be an integer coprime to $p$ and
let $X$ be either the affine modular curve $Y_1(N)$ over $R$ or its ordinary locus
$Y_1(N)_{ord}$ (i.e. the locus where the Eisentein form $E_{p-1}$ is invertible). Let $L$ be the line bundle on the complete modular curve $X_1(N)$ over $R$
such that the global sections of the powers  $L^{\otimes \k}$, $\k \geq 0$,  are the classical modular forms (on $\Gamma_1(N)$) of  weight $\k$ over $R$ and let
$V\ra X$, $V:=Spec\ \bigoplus_{\k \in {\mathbb Z}}L^{\otimes  \k}$, be the natural ${\mathbb G}_m$-torsor associated to the restriction of $L$ to $X$. A {\it $\d$-modular function of order $r$} (on $\Gamma_1(N)$)  \cite{difmod,book} will mean a $\d$-function of order $r$ on $V$, i.e. an element of $\cO^r(V)$. Let $W={\mathbb Z}[\phi]$ be the polynomial ring in the variable $\phi$. Then the multiplicative monoid $W$ naturally acts on $R^{\times}$; for $w \in W$ and $\lambda\in R^{\times}$ we write
$(w,\lambda)\mapsto \lambda^w$ for the action. Evaluation at $\phi=1$ defines a ring homomorphism $deg:W={\mathbb Z}[\phi]\ra {\mathbb Z}$.
A $\d$-modular function $f\in \cO^r(V)$ will be called a {\it $\d$-modular form of weight} $w \in W$ if
$f(\lambda \cdot a)=\lambda^w f(a)$ for $a\in V(R)$ and $\lambda \in R^{\times}$, where $\lambda \cdot a$ is defined via the ${\mathbb G}_m$-action on $V$.

\subsection{$\d$-Fourier expansions}
Any $\d$-modular function of order $r$ has a natural {\it $\d$-Fourier expansion} in the ring
of {\it $\d$-series} $R((q))[q',...,q^{(r)}]\h$
where $q,q',...,q^{(r)}$ are variables, $R((q)):=R[[q]][1/q]$, and the upper $\h$ means here (and everywhere later) {\it completion in the $p$-adic topology}.

There are unique $p$-derivations $\d$ from $R((q))[q',...,q^{(r)}]\h$
to $R((q))[q',...,q^{(r+1)}]\h$ extending $\d$ on $R$ and such that $\d q=q'$, $\d q'=q''$, etc.
The $\d$-Fourier expansion maps are compatible with the classical Fourier expansion maps and commute with $\d$. Recall that for $\k \in {\mathbb Z}_{\geq 0}$ the classical Hecke operators $T_{\k+2}(n)$ (with $n \geq 1$, $(n,p)=1$) and $T_{\k+2}(p)$ act on $R((q))$. We have an induced action of $T_{\k+2}(n)$, $T_{\k+2}(p)$ on $k((q))$; clearly $T_{\k+2}(p)$ on $k((q))$
coincides with Atkin's operator $U$ on $k((q))$, defined by $U(\sum a_n q^n)=\sum a_{np}q^n$.
A series $\varphi \in k((q))$ will called {\it primitive} if $U \varphi=0$. A $\d$-series in $k((q))[q',...,q^{(r)}]$ will be called {\it primitive}
if its image in $k((q))$ under the specialization $q'=...=q^{(r)}=0$ is primitive.
One can define Hecke operators $T_{\k}(n)$, $pT_{\k}(p)$ on $R((q))[q',...,q^{(r)}]\h$ (where $pT_{\k}(p)$ is only ``partially defined" i.e. defined on an appropriate subspace); cf. Sections 2 and 3 below for all the relevant details. These operators induce operators $T_{\k}(n)$, ``$pT_{\k}(p)$" on $k((q))[q',...,q^{(r)}]$ (where ``$pT_{\k}(p)$" is only ``partially defined" i.e. defined on an appropriate subspace; the ``\ " signs are meant to remind us that the operator $T_{\k}(p)$ itself is not defined mod $p$).

\subsection{Main result}
The following is our main result; it is a consequence of Theorems \ref{dezastru} and \ref{converseofdezastru} in the body of the paper.
Assume $X=Y_1(N)_{ord}$ and let $\k \in {\mathbb Z}_{\geq 0}$.

\begin{theorem}
\label{maint}
 There is a one-to-one correspondence between the following sets of objects:

i) Series in $qk[[q]]$  which are eigenvectors of all  Hecke operators $T_{\k+2}(n)$, $T_{\k+2}(p)$, $(n,p)=1$, and which are Fourier expansions of classical modular forms  over $k$ of weight $\equiv \k+2$ mod $p-1$;

ii) Primitive  $\d$-series in $k[[q]][q']$ which are eigenvectors of all Hecke operators $nT_{\k}(n)$, ``$pT_{\k}(p)$", $(n,p)=1$, and which are $\d$-Fourier expansions of $\d$-modular forms of some order $r \geq 0$ and  weight $w$ with $deg(w)=\k$.

This correspondence preserves the respective eigenvalues.
\end{theorem}

\begin{remark}

1) As Theorems \ref{dezastru} and \ref{converseofdezastru} will show the correspondence in Theorem \ref{maint} is given, on a computational level, by an entirely explicit formula (but note that the proof that this formula establishes the desired correspondence is {\it not} merely computational.) The formula is as follows. If $\varphi=\sum_{m \geq 1} a_m q^m \in k[[q]]$ is a series as in i) of the Theorem then $a_1 \neq 0$ and the corresponding $\d$-series in ii) is given by
$$\varphi^{\sharp,2}:=\sum_{(n,p)=1}\frac{a_n}{n} q^n -\frac{a_p}{a_1} \cdot \left(\sum_{m\geq 1}a_m q^{mp}\right) \frac{q'}{q^p}+ e \cdot \left(\sum_{m\geq 1}a_m q^{mp^2}\right) \cdot \left(\frac{q'}{q^p}\right)^p\in k[[q]][q'],$$
where $e$ is $1$ or $0$  according as $\k$ is $0$ or $>0$. (The upper index $2$ in $\varphi^{\sharp,2}$ is meant to reflect the $p^2$ exponent in the right hand side of the above equality; later in the body of the paper we will encounter a $\varphi^{\sharp,1}$ series as well. The $\sharp$ sign is meant to reflect the link between these objects and the objects $f^{\sharp}$ introduced in \cite{eigen}.)

2)  Theorem \ref{maint} provides a  complete description
of primitive  $\d$-series mod $p$ of order $1$ which are eigenvectors of all the Hecke operators and which are $\d$-Fourier expansions of $\d$-modular forms of arbitrary order. It would be desirable to have such a description in characteristic zero
and/or for higher order $\d$-series.
However note that all known examples (so far) of $\d$-modular forms of order $\geq 2$ which are eigenvectors of all Hecke operators  have the property that their $\d$-Fourier expansion reduced mod $p$ has order $1$; by the way some of these forms play a key role in \cite{eigen,dcc,BP}. So it is reasonable to ask if it is true that  {\it any $\d$-modular form of order $\geq 1$ which is an eigenvector of all the Hecke operators must have a $\d$-Fourier expansion whose reduction mod $p$ has order $1$}.

3) Note that in ii) of the above Theorem one can take the order to be $r=1$ and the weight to be
$w=\k$. Also note that the $\d$-modular forms in ii) above have, a priori,
``singularities" at the cusps and at the supersingular points. Nevertheless, in the special case when the classical modular forms in i) above come from newforms on $\Gamma_0(N)$ over ${\mathbb Z}$ of weight $2$   one can choose the $\d$-modular forms in ii) of weight $0$, order $2$, and {\it without singularities at the cusps or at at the supersingular points}; this was done in \cite{eigen} where the corresponding $\d$-modular forms were denoted
(at least in the ``non-CL" case) by $f^{\sharp}$. These $f^{\sharp}$s played, by the way, a key role in the proof of the main results in \cite{BP} about linear dependence relations among Heegner points.
It would be interesting to find analogues of the forms $f^{\sharp}$ in higher weights.

4) One of the subtleties of the above theory is related to the fact that the operator ``$pT_{\k}(p)$" is not everywhere defined.
The failure of this operator to be everywhere defined is related to the failure of ``the fundamental theorem of symmetric polynomials"
 in the context of $\d$-functions; cf. \cite{eigen,dcc}. The domain of definition of ``$pT_{\k}(p)$" will be the space of all $\d$-series
for which the analogue of ``the fundamental theorem of symmetric polynomials" holds; these $\d$-series will be called {\it Taylor $\d-p$-symmetric}. One of our main results  will be a complete determination the space of Taylor $\d -p$-symmetric $\d$-series; cf. Theorems \ref{conv} and \ref{infinitesums}.

\end{remark}

\subsection{Comments on $\d$-geometry \cite{book}}
The present paper fits into a more general program for which we refer to \cite{book}. Roughly speaking this program proposes to enrich (usual) algebraic geometry by replacing algebraic equations (i.e. expressions of the form $f=0$, $f$ a polynomial function) with arithmetic differential equations
(i.e. expressions of the form $f=0$, $f$ a $\d$-function). This enriched geometry can be referred to as {\it $\d$-geometry}. One of the main motivations/applications of $\d$-geometry is the construction of  certain quotients of (usual) algebraic curves by actions of (usual) correspondences. Such quotients fail to exist within (usual) algebraic geometry in the sense that the corresponding categorical quotients in (usual) algebraic geometry reduce to a point. On the contrary, in $\d$-geometry, one can construct a number of interesting such categorical quotients, e.g. the quotient of the modular curve $Y_1(N)$ by the action of the Hecke correspondences. The construction/underdstanding of the latter is based upon the theory of $\d$-modular forms.

On a more ``philosophical" level note that $\d$-geometry and, more generally, $\Lambda$-geometry (which is a  several prime generalization of $\d$-geometry)  can be viewed as an incarnation of  the  ``geometry over the field with one element";
cf. the Introduction to \cite{book} for  remarks on the  single prime case and \cite{f1} for a systematic explanation of this viewpoint  in  the several prime case.

On the other hand, from a more ``pragmatic" point of view, we note that  $\d$-geometry has  applications to (usual) arithmetic geometry such as: matters surrounding the  Manin-Mumford conjecture \cite{pjets,dcc}, congruences between (usual) modular forms
\cite{difmod,Barcau}, and linear dependence relations among Heegner points \cite{BP}.

\subsection{Plan of the paper} Sections 2 and 3   introduce Hecke operators $T_{\k}(n)$, $(n,p)=1$ and ``$pT_{\k}(p)$"
respectively, acting on $\d$-series. Section 4 gives the complete determination of the $\d$-series mod $p$ of order $1$ for which ``the analogue of the fundamental theorem of symmetric polynomials" holds. Section 5 gives a multiplicity one theorem for $\d$-series which are eigenvectors of all Hecke operators.
Section 6 begins with an overview of $\d$-modular forms \cite{difmod, book}  and Serre-Katz   $p$-adic modular forms \cite{Katz};
 then we use the multiplicity one result plus  results in \cite{difmod,book} and \cite{Katz} to prove results implying Theorem \ref{maint}.

\subsection{Acknowledgment} While writing this paper the first author was
partially supported by NSF grant  DMS-0852591.  Any opinions,
findings, and conclusions or recommendations expressed in this
material are those of the authors and do not necessarily reflect the
views of the National Science Foundation.

\section{Hecke operators away from $p$}

\subsection{Classical Hecke operators}
Throughout the paper the divisors  of a given non-zero integer  are always taken to be positive, the greatest common divisor of two non-zero integers $m,n$ is denoted by $(m,n)$, and we use the convention $(m,n)=n$ for $m=0$, $n \neq 0$. Fix throughout the paper an integer  $N \geq 4$  and let $\epsilon:{\mathbb Z}_{> 0}\ra \{0,1\}$ be the ``trivial primitive character" mod $N$ defined by $\epsilon(A)=1$ if $(A,N)=1$ and $\epsilon(A)=0$ otherwise.

For each integer $n \geq 1$  and each integer $N \geq 4$
 consider the set
$$\{(A,B,D); A,B,D \in {\mathbb Z}_{\geq 0}, AD=n, (A,N)=1, B<D\}$$
Triples $A,B,D$ will always be assumed to be in the set above. Recall (cf., say, \cite{Knapp})
the action of the $n$-th Hecke operator $T_{\k}(n)$ on classical modular forms $f=\sum_{m \geq 0}a_mq^m$ on $\Gamma_0(N)$ of weight $\k\geq 2$ with complex coefficients $a_m \in {\mathbb C}$ given by

\bigskip

$$\begin{array}{rcl}
T_{\k}(n)f & := & n^{\k-1}\sum_{A,B,D}D^{-\k} f(\zeta_D^B q^{A/D})\\
\  & \  & \  \\
& = & \sum_{m \geq 0}\left(
\sum_{A|(n,m)} \epsilon(A)A^{\k-1}a_{\frac{mn}{A^2}}\right) q^m.
\end{array}$$
\bigskip

\noindent Here $q=e^{2 \pi \sqrt{-1}z}$, $\zeta_D:=e^{2 \pi \sqrt{-1}/D}$.

\subsection{Hecke operators $T_{\k}(n)$ on $\d$-series}
Now assume $n$ and $N$ are coprime to $p$ and assume $q,q',q'',...,q^{(r)},...$ are indeterminates.

\begin{definition}
For each integer $\k \in {\mathbb Z}$  the {\it Hecke} operator $f \mapsto T_{\k}(n)f$ on $R((q))[q',...,q^{(r)}]\h$ is defined as follows.
For $f=f(q,q',...,q^{(r)})$,
\begin{equation}
\label{eqq}
T_{\k}(n)f:=n^{\k-1}\sum_{A,B,D}D^{-\k} f(\zeta_D^B q^{A/D}, \d(\zeta_D^B q^{A/D}),...,\d^r(\zeta_D^B q^{A/D})).\end{equation}
\end{definition}

Here $\zeta_D=\zeta_n^{n/D}\in R$ where $\zeta_n \in R$ is a fixed primitive $n$-th root of unity and the right hand side of (\ref{eqq})
is a priori in the ring
\begin{equation}
\label{rring}
R((q_n))\h[q'_n,...,q_n^{(r)}]\h,\ \ q_n=q^{1/n}.\end{equation}
However, by \cite{book} Proposition 3.13,
$$q_n',...,q_n^{(r)} \in R[q,q^{-1},q',...,q^{(r)}]\h$$
hence the ring (\ref{rring}) equals
$$R((q_n))\h[q',...,q^{(r)}]\h.$$
Since $T_{\k}(n)f$ is invariant under the substitution $q_n^{(i)}\mapsto \d^i(\zeta_n q_n)$ it follows that
$T_{\k}(n)f \in R((q))\h[q',...,q^{(r)}]\h$. So the operators $T_{\k}(n)$ send
$R((q))\h[q',...,q^{(r)}]\h$ into itself.
As we shall see below for $n \geq 2$ the operators $T_{\k}(n)$ do {\it not} send $R[[q]][q',...,q^{(r)}]\h$ into itself.

The operators $T_{\k}(n)$ on $R((q))[q',...,q^{(r)}]\h$ induce operators still denoted by $T_{\k}(n)$
on $k((q))[q',...,q^{(r)}]$.

Recall the operator $V$ on $R((q))\h$ defined by $V(\sum a_nq^n)=\sum a_nq^{pn}$.
It induces an operator still denoted by $V$ on $k((q))$.

For $r=0$, $T_{\k}(n)$  commute with the operator $V$ on $R((q))\h$.

\subsection{Order $r=1$}
We have the following formula for the Hecke action on $\d$-series of order $1$:

\begin{proposition}
\label{Tf}
Assume that
\begin{equation}
\label{f}
f=\sum_{m,m'} a_{m,m'} q^m(q')^{m'}\end{equation}
where $m \in \bZ$, $m'\in \bZ_{\geq 0}$. Then we have the following congruence mod $(p)$:
\begin{equation}
 \label{cuani}
 T_{\k}(n)f \equiv \sum_{m,m'} \left( \sum_{A | (n,m)} n^{-m'}\epsilon(A)A^{\k+2m'-1}
a_{\frac{mn}{A^2}-m'p,m'} \right) q^{m-m'p}(q')^{m'}.\end{equation}
\end{proposition}

{\it Proof}. Note that
\begin{equation}
\label{ord1}
\begin{array}{rcl}
\d(\zeta_D^B q^{A/D}) & = & \frac{1}{p}[\phi(\zeta_D^B q^{A/D})-(\zeta_D^B q^{A/D})^p]\\
\  & \  & \ \\
\  & = & \frac{1}{p} [\zeta_D^{Bp} (q^p+pq')^{A/D}-\zeta_D^{Bp} q^{Ap/D}]\\
\  & \  & \ \\
\  & \equiv & \frac{A}{D}\zeta_D^{Bp} q^{(A-D)p/D} q'\ \ \ mod\  (p).\end{array}\end{equation}
Then the formula in the statement of the Proposition follows by a simple computation, using  the fact that
$$\sum_{B=0}^{D-1} \zeta_D^{m+m'p}$$
is $D$ or $0$ according as $D$ divides or does not divide $m+m'p$.
\qed

\bigskip

\begin{corollary}
\label{buium}
Let
\begin{equation}
\label{theseries}
 \overline{f}=\sum_{m'} \overline{f}_{m'}(q) \left( \frac{q'}{q^p} \right)^{m'}\in k((q))[q'],\ \ \ \overline{f}_{m'}(q)\in k((q)).\end{equation}
Then for any integer $\k$ and any integer  $n \geq 1$  coprime to $p$ we have:
$$T_{\k}(n)\overline{f}= \sum_{m'} n^{-m'}(T_{\k+2m'}(n)\overline{f}_{m'}(q))\left( \frac{q'}{q^p} \right)^{m'}.$$
In particular for
  $\overline{\lambda}_n \in k$ we have  $T_{\k}(n)\overline{f}= \overline{\lambda}_n \overline{f}$ if and only if
$$T_{\k+2m'} (n) \overline{f}_{m'}= n^{m'}\overline{\lambda}_n \overline{f}_{m'}\ \   \text{for  all}\ \  m'\geq 0.$$
\end{corollary}

{\it Proof}. This follows immediately from Proposition \ref{Tf}.\qed

Let us say that a series in $k((q))[q',...,q^{(r)}]$  is {\it holomorphic  at infinity} if it belongs to $k[[q]][q',...,q^{(r)}]$. Also denote by $v_p$ the $p$-adic valuation on ${\mathbb Z}$.

\begin{corollary}
\label{echelon}
Assume that, for a given $\k\in {\mathbb Z}$ the series $\overline{f}\in k[[q]][q']$   has the property that $T_{\k}(n)\overline{f}$ is holomorphic at infinity for all $n \geq 1$ coprime to $p$.
Then $\overline{f}$ has the form
\begin{equation}
\label{wash}
\overline{f}(q,q')  =  \varphi_0(q)+\sum_{m'\geq 1} (V^{v_p(m')+1} (\varphi_{m'}(q)))\left( \frac{q'}{q^p}\right)^{m'},\end{equation}
with
\begin{equation}
\label{hair}
\varphi_0  \in  k[[q]],\ \ \
\varphi_{m'}(q)  \in  q^{m'/p^{v_p(m')}}k[[q]]\ \ \text{for}\ \ m'\geq 1.\end{equation}
\end{corollary}

{\it Proof}. Note that, since $T_{\k}(1)\overline{f}=\overline{f}$, $\overline{f}$ is holomorphic at infinity
so equation (\ref{hair}) follows from (\ref{wash}).  Let $\overline{f}$ be the reduction mod $p$ of a series as in (\ref{f}).
It is enough  to show if two integers $m_0 \geq 1$ and $m'\geq 1$
satisfy $v_p(m_0)\leq v_p(m')$
then $\overline{a}_{m_0,m'}=0$. Pick such integers $m_0,m'$ and set $i=v_p(m_0)$, $m_0=p^i\mu$, $m'=p^i \mu'$, $n=\mu+p\mu'$. Clearly $n$ is coprime to $p$.
 Picking out the coefficient of $q^{p^i-p^{i+1}\mu'}(q')^{p^i\mu'}$ in the equation in
 Proposition \ref{Tf} we get
$$\overline{a}_{m_0,m'}=\overline{a}_{p^i n-p^{i+1}\mu',p^i \mu'}= 0$$
 and we are done. \qed

\bigskip

\begin{corollary}
\label{voce}
Let $\k$ be an integer, let $\overline{f}\in k[[q]][q']$  be holomorphic  at infinity, and assume that
for any integer $n \geq 1$ coprime to $p$ we are given a $\overline{\lambda}_n \in k$. Then  $T_{\k}(n)\overline{f}= \overline{\lambda}_n \overline{f}$ for all
$(n,p)=1$ if and only if
 $\overline{f}$ has the form (\ref{wash}) and
 $$T_{\k+2m'} (n) \varphi_{m'}(q)= n^{m'}\overline{\lambda}_n \varphi_{m'}(q)\ \ \text{for all}\ \ m'\geq 0.$$
\end{corollary}

{\it Proof}. This follows directly from the previous corollaries plus the commutation of $T_{\k}(n)$ and $V$ on $k[[q]]$.
\qed

\subsection{Order $r=2$}
Let us record the  formula giving the Hecke action on  $\d$-series of order $2$.
This formula will not be used in the sequel.

\begin{proposition}
\label{ord22}
If $f=\sum_{m,m
,m''} a_{m,m',m''} q^m (q')^{m'}(q'')^{m''}\in R((q))[q',q'']\h$ then we have the following congruence mod $p$:
$$\begin{array}{rcl}
T_{\k}(n)f & \equiv & \sum
A^{\k-1}  \left(\frac{A}{D}\right)^{m'+m''} \times a_{m,m',m''} \times q^{A(m+m'p+m''p^2)/D} \\
\  & \  & \  \\
\  & \  &  \times  \left(\frac{q'}{q^p}\right)^{m'} \times \left[ \frac{q''}{q^{p^2}}+\frac{\d(A/D)}{A/D} \cdot \left(\frac{q'}{q^p}\right)^p+
\frac{1}{2}\left( \frac{A}{D}-1\right)  \cdot \left( \frac{q'}{q^p}\right)^{2p}\right]^{m''}\end{array}$$
where the sum in the right hand side  runs through all $m,m',m'',A,D$ with $A \geq 1, AD=n, (A,N)=1,D|m+m'p+m''p^2$.
\end{proposition}

{\it Proof}.
A computation similar to the one in the proof of Proposition \ref{Tf}.
\qed

\medskip

Note that  the formula in Proposition \ref{ord22} acquires a simpler form for special $n$s. Indeed assume  $n=\ell$ is a prime. If $\ell \equiv 1$ mod $p$ then
$\frac{A}{D}-1=0$ in $k$. If $\ell \equiv 1$ mod $p^2$ then $\d(A/D)=0$ in $k$.
Finally if $\ell \equiv 1$ mod $p$ but $\ell \not\equiv 1$ mod $p^2$ then $\d(A/D)\neq 0$ in $k$.

\subsection{Frobenii}
\label{Frobenii}
Consider the ring endomorphisms $F,F_k,F_{/k}$  of $k((q))[q',...,q^{(r)}]$ defined as follows:
$F$ is the $p$-power Frobenius (the ``absolute Frobenius"); $F_k$ is the ring automorphism that acts as the $p$-power Frobenius on $k$ and is the identity on the variables $q,q',...,q^{(r)}$; $F_{/k}$ is the ring endomorphism that is the identity on $k$ and sends $q,q',...,q^{(r)}$ into $q^p,(q')^p,...,(q^{(r)})^p$ respectively (the ``relative Frobenius"). So we have $F=F_k \circ F_{/k}=F_{/k}\circ F_k$. Of course
$V=F_{/k}$ on $k((q))$. Also clearly $T_{\k}(n)$ commute with $F$. By Proposition \ref{Tf} $T_{\k}(n)$ also commute with  $F_{k}$ on $k((q))[q']$; so
$T_{\k}(n)$  commute with $F_{/k}$ on $k((q))[q']$.

\bigskip

\section{Hecke operator at $p$}

\subsection{Taylor and Laurent $\d$-symmetry}
Following  \cite{eigen} we consider   the $R-$algebras
$$
\begin{array}{rcl}
A & := &
R[[s_1,...,s_p]][s_p^{-1}]\h[s_1',...,s_p',...,s_1^{(r)},...,s_p^{(r)}]\h,\\
\  & \  & \  \\
B & := &
R[[q_1,...,q_p]][q_1^{-1}...q_p^{-1}]\h[q_1',...,q_p',...,q_1^{(r)},...,q_p^{(r)}]\h,
\end{array}$$
where $s_1,...,s_p,s'_1,...,s_p',...$ and $q_1,...,q_p,q'_1,...,q_p',...$ are indeterminates. In \cite{eigen}, Lemma 9.10
 we proved that the natural algebra map $$A \ra B,\ \ \ s^{(i)}_j
\mapsto \d^iS_j,$$
where   $S_1,...,S_p$ are the fundamental symmetric polynomials in $q_1,...,q_p$,
 is injective with torsion free cokernel.
 We will view this algebra map as an inclusion.

 \begin{definition}
 An element  $G \in B$ is
 called {\it Laurent $\d-$symmetric} \cite{eigen}
 if it is the image of some
element $G_{(p)} \in A$ (which is then unique). An element $f \in R((q))\h[q',...,q^{(r)}]\h$ will be called {\it Laurent $\d-p$-symmetric} if
\[\Sigma_pf:=\sum_{j=1}^pf(q_j,...,q_j^{(r)}) \in B\]
is Laurent $\d-$symmetric.
\end{definition}

In the same way one can consider the algebras
$$
\begin{array}{rcl}
A & := &
R[[s_1,...,s_p]][s_1',...,s_p',...,s_1^{(r)},...,s_p^{(r)}]\h,\\
\  & \  & \  \\
B & := &
R[[q_1,...,q_p]][q_1',...,q_p',...,q_1^{(r)},...,q_p^{(r)}]\h.
\end{array}$$ As before the natural algebra map $$A \ra B,\ \ \ s^{(i)}_j
\mapsto \d^iS_j,$$
 is injective with torsion free cokernel.

 \begin{definition}
 An element  $G \in B$ will be
 called {\it  Taylor $\d-$symmetric} if it is the image of some
element $G_{(p)} \in A$ (which is then unique). An element $f \in R[[q]][q',...,q^{(r)}]\h$ will be called {\it Taylor $\d-p$-symmetric} if
\[\Sigma_pf:=\sum_{j=1}^pf(q_j,...,q_j^{(r)}) \in B\]
is  Taylor $\d-$symmetric.
\end{definition}

Clearly a Taylor $\d-p$-symmetric series is also Laurent $\d-p$-symmetric.

\begin{remark}

1) Any element of $R[[q]]$ (respectively $R((q))$) is Taylor (respectively Laurent) $\d-p$-symmetric.

2) The Taylor (respectively Laurent) $\d-p$-symmetric   elements in $R[[q]][q',...,q^{(r)}]\h$
(respectively $R((q))\h[q',...,q^{(r)}]\h$) form a $p$-adically closed $R$-submodule.

3) If $f$ is Taylor (respectively Laurent) $\d-p$-symmetric  then $\phi(f)$ is Taylor (respectively Laurent) $\d-p$-symmetric.

4) If $f\in R[[q]][q',...,q^{(r)}]\h$ (respectively $f \in R((q))\h[q',...,q^{(r)}]\h$) and $pf$ is
Taylor (respectively Laurent) $\d-p$-symmetric
 then $f$ is Taylor (respectively Laurent) $\d-p$-symmetric.

5) By 1)-4) any element $f$ in $R[[q]][q',...,q^{(r)}]\h$ (respectively in $R((q))\h[q',...,q^{(r)}]\h$) of the form
$$f=\frac{\sum_{i=0}^m \phi^i(g_i)}{p^{\nu}}$$
where $g_i$ are in $R[[q]]$ (respectively in $R((q))$)  is
 Taylor (respectively Laurent) $\d-p$-symmetric.
 In particular for any $g$ in $R[[q]]$ (respectively in $R((q))$) we have that
$\d g=\frac{\phi(g)-g^p}{p}$, and more generally $\frac{\phi^i(g)-g^{p^i}}{p}$ are Taylor (respectively Laurent) $\d-p$-symmetric.

6) Let  $\cF \in R[[T_1,T_2]]^g$ be a formal group law,
and let $\psi \in R[[T]][T,...,T^{(r)}]\h$ be such that
\[\psi(\cF(T_1,T_2),...,\d^r\cF(T_1,T_2))=
\psi(T_1,...,T_1^{(r)})+\psi(T_2,...,T_2^{(r)})\] in the
ring
\[R[[T_1,T_2]][T'_1,T_2',,...,T_1^{(r)},T_2^{(r)}]\h.\]
(Such a $\psi$ is called a $\d$-{\it character} of ${\mathcal F}$.)
Let $\varphi(q) \in qR[[q]]$ and let
\[f:=\psi(\varphi(q),...,\d^r(\varphi(q))) \in
R[[q]][q',...,q^{(r)}]\h.\] Then $f$ is Taylor $\d-p$-symmetric. Cf the argument in \cite{dcc}.

Note that if $\cF$ is defined over ${\mathbb Z}_p$  then $\cF$ posses a  $\d$-character $\psi$ of order $r$
at most the height of $\cF$ mod $p$ such that
$$\psi(T,0,...,0)\in T+T^p{\mathbb Z}_p[[T]];$$
 cf. \cite{book}, proof of Proposition 4.26.

 Applying the above considerations to the multiplicative formal group we get that for any $\varphi(q)\in qR((q))$ the series
 $$\frac{1}{p}\log\left(\frac{\phi(\varphi(q)+1)}{(\varphi(q)+1)^p}\right)$$
 is Taylor $\d-p$-symmetric. (Here, as usual,  $\log(1+T)=T-T^2/2+T^3/3-...$)

 7) The series
 \begin{equation}
 \label{Psi}\Psi=\frac{1}{p}\log \left(\frac{\phi(q)}{q^p}\right)\end{equation}
 is Laurent $\d-p$-symmetric; cf. \cite{eigen},  proof of Proposition 9.13.

 8) In \cite{eigen} we also defined the concept of $\d$-symmetric element in
 $$R[[q_1,...,q_p,...,q_1^{(r)},...,q_p^{(p)}]]$$
  (without the qualification ``Taylor" or ``Laurent"). We will not use this concept in the present paper.
  But note that if a series is Taylor $\d$-symmetric then it is also $\d$-symmetric in the sense of \cite{eigen} (and Laurent $\d$-symmetric in the sense of the present paper).
\end{remark}

\begin{definition}
For any  Taylor (respectively Laurent) $\d-p$-symmetric  $$f\in R[[q]][q',...,q^{(r)}]\h \ \ \text{(respectively
  $f \in R((q))\h[q',...,q^{(r)}]\h$)}$$ we  define
$$
 U f:=p^{-1}(\Sigma_p f)_{(p)}(0,...,0,q,...,0,...,0,q^{(r)})$$
which is  an element in
$p^{-1}R[[q]][q',...,q^{(r)}]\h$ (respectively
 in $p^{-1}R((q))\h[q',...,q^{(r)}]\h$). \end{definition}

 The operator $pU$ takes
 $R[[q]][q',...,q^{(r)}]\h$ (respectively
 in $R((q))\h[q',...,q^{(r)}]\h$) into $R[[q]][q',...,q^{(r)}]\h$ (respectively
 in $R((q))\h[q',...,q^{(r)}]\h$).
 On the other hand the  restriction of $U$ to $R((q))\h$ (respectively $R[[q]]$) takes values in $R((q))\h$ (respectively $R[[q]]$) and is equal to the classical $U$-operator
$$U(\sum a_mq^m)=\sum a_{mp}q^m.$$

\begin{definition}
Define for any $f \in R((q))\h[q',...,q^{(r)}]\h$ the series
$$
Vf  := f(q^p,...,\d^r(q^p)) \in R((q))\h[q',...,q^{(r)}]\h.$$
So for any Taylor (respectively Laurent) $\d-p$-symmetric  $f$ in $R[[q]][q',...,q^{(r)}]\h$ (respectively
 in $R((q))\h[q',...,q^{(r)}]\h$) and any $\k\in {\mathbb Z}$ we may define
$$pT_{\k}(p)f  =  pUf+p^{\k}Vf$$
which is  an element in $p^{\k}R[[q]][q',...,q^{(r)}]\h$ (respectively
 in $p^{\k}R((q))\h[q',...,q^{(r)}]\h$). \end{definition}

 The restriction of $pT_{\k}(p)$ to $R((q))$
 is, of course, $p$ times the ``classical" Hecke operator $T_{\k}(p)$ on $R((q))$ defined by
 $$T_{\k}(p)(\sum a_m q^m)=\sum a_{pm}q^m+p^{\k-1}\sum a_m q^{pm}.$$
Recall:

\begin{proposition}
\cite{eigen}
The series $\Psi$ in (\ref{Psi}) satisfies
$$pU\Psi=\Psi,\ \ V\Psi=p\Psi.$$
\end{proposition}

For the next definition recall that the homomorphism
$$\overline{A}:=A\otimes_R k \ra \overline{B}:=B \otimes_R k$$
is injective (in both situations described in the beginning of the section).

\begin{definition}
An element $\overline{G} \in \overline{B}$ is
 called {\it Taylor $\d$-symmetric mod $p$} (respectively {\it Laurent $\d-$symmetric mod $p$})
  if it is the image of some
element $\overline{G}_{(p)} \in \overline{A}$ (which is then unique). An element $\overline{f}\in k[[q]][q',...,q^{(r)}]\h$ (respectively $\overline{f} \in k((q))[q',...,q^{(r)}]$) will be called Taylor (respectively Laurent) $\d-p$-{\it symmetric}  if
\[\Sigma_p\overline{f}:=\sum_{j=1}^p\overline{f}(q_j,...,q_j^{(r)}) \in \overline{B}\]
is  Taylor $\d-$symmetric mod $p$ (respectively Laurent $\d$-symmetric mod $p$).
\end{definition}

Clearly any Taylor $\d-p$-symmetric series is Laurent $\d-p$-symmetric.

\begin{remark}
\label{multe}
1) The Taylor (respectively Laurent) $\d-p$-symmetric   elements in $k[[q]][q',...,q^{(r)}]$ (respectively in $k((q))[q',...,q^{(r)}]$)
form a $k$-subspace closed under $F_k$ and $F$ (hence also under $F_{/k}$).

2) If $f \in R[[q]][q',...,q^{(r)}]\h$ (respectively $f\in R((q))\h[q',...,q^{(r)}]\h$) is congruent mod $p$  to a
Taylor (respectively Laurent) $\d-p$-symmetric   element then the image of $\overline{f}$ of $f$ in $k[[q]][q',...,q^{(r)}]$ (respectively in $k((q))[q',...,q^{(r)}]$)
Taylor (respectively Laurent) $\d-p$-symmetric.\end{remark}

\begin{definition}
For any Taylor (respectively Laurent) $\d-p$-symmetric
$$\text{$\overline{f}\in k[[q]][q',...,q^{(r)}]\h$ (respectively $k((q))[q',...,q^{(r)}]$)}$$
we may define
$$
``pU" \overline{f}:=(\Sigma_p \overline{f})_{(p)}(0,...,0,q,...,0,...,0,q^{(r)})$$
which is an element of $k[[q]][q',...,q^{(r)}]\h$ (respectively $k((q))[q',...,q^{(r)}]$).
\end{definition}

The operator $``pU"$ clearly commutes with the operators $F$ and $F_k$ and hence it also commutes with the operator $F_{/k}$ (cf. section \ref{Frobenii}).
If
$$\text{$f \in R[[q]][q',...,q^{(r)}]\h$ (respectively $f \in R((q))\h[q',...,q^{(r)}]\h$)}$$
 is Taylor
(respectively Laurent) $\d-p$-symmetric and $\overline{f}$ is the reduction mod $p$  of $f$ viewed as an element in $k[[q]][q',...,q^{(r)}]$ (respectively in $k((q))[q',...,q^{(r)}]$) then $``pU"\overline{f}$ is the reduction mod $p$ of $pUf$; this justifies the notation in $``pU"\overline{f}$.

Note that the operator $U:R((q))\h\ra R((q))\h$ induces an operator still denoted by $U$,
$U:k((q))\ra k((q))$ (which is, of course, the classical $U$-operator $U\overline{f}=
\sum \overline{a}_{mp}q^m$, for $\overline{f}=\sum \overline{a}_m q^m\in k((q))$). On the other hand note that
$``pU"\overline{f}=0$ for all $\overline{f}\in k((q))$. Finally note that if $\k \geq 1$ then the operator $T_{\k}(p)$ on $R((q))$
induces an operator $T_{\k}(p)$ on $k((q))$; if $\k \geq 2$ then $T_{\k}(p)$ on $k((q))$ coincides with $U$ on $k((q))$.

\begin{definition}
Define   the ring endomorphism $V$ of
$$
\text{$k[[q]][q',...,q^{(r)}]$ (respectively $k((q))[q',...,q^{(r)}]$)}$$ as the reduction mod $p$ of the operator
 $V$ over $R$. (Note that $V(q')=0$ and $F_{/k}(q')=(q')^p$ so in particular $V\neq F_{/k}$ on $k((q))[q']$.)
As in the case of characteristic zero, for any $\kappa\in {\mathbb Z}_{\geq 0}$ and any Taylor (respectively Laurent) $\d-p$-symmetric series $\overline{f}$ in $k[[q]][q',...,q^{(r)}]$ (respectively $k((q))[q',...,q^{(r)}]$) we define
$$``pT_{\k}(p)"\overline{f}=``pU"\overline{f}+\ee \cdot V\overline{f}$$
which is again an element of $k[[q]][q',...,q^{(r)}]$ (respectively $k((q))[q',...,q^{(r)}]$).
(Note that $\ee$ is $0$ or $1$ according as $\k$ is $>0$ or $0$.)
 \end{definition}

 The operator $V$ clearly commutes with $F$ and $F_k$ (and hence also with $F_{/k}$).
 So the operators $``pT_{\k}(p)"$ commute with $F,F_k,F_{/k}$.

 Also for $f$ any Taylor (respectively Laurent) $\d-p$-symmetric series  in $R[[q]][q',...,q^{(r)}]\h$ (respectively $R((q))\h[q',...,q^{(r)}]\h$) with reduction mod $p$ $\overline{f}$ we have that $``pT_{\k}(p)"\overline{f}$ is the reduction mod $p$ of $pT_{\k}(p)f$ which, again, justifies our notation.

\section{Structure of Laurent and Taylor $\d-p$-symmetric series}
 In what follows we address the problem of determining what series are Laurent (respectively Taylor)
 $\d-p$-symmetric and determining the action of our operators $``pU"$ on them. We will use the following notation:
 for all $\varphi=\sum \overline{a}_nq^n \in k((q))$ we define
 \begin{equation}
 \label{varphiminunu}\varphi^{(-1)}:=\theta^{p-2}\varphi=\sum_{(n,p)=1} \frac{\overline{a}_n}{n}q^n \in k((q))\end{equation}
 where $\theta=q\frac{d}{dq}$ is the Serre theta operator.

 \begin{theorem}
\label{conv}
If an element $\overline{f}\in k[[q]][q']$
is Taylor $\d-p$-symmetric then  it has the form
\begin{equation}
\label{blabla}
\overline{f}=\varphi_0(q)+\sum_{s\geq 0} (V^{s+1}(\varphi_{p^s}(q)))\left( \frac{q'}{q^p}\right)^{p^s}\in k((q))[q']\end{equation}
with $\varphi_0(q)\in k[[q]]$, $\varphi_1(q),\varphi_p(q),\varphi_{p^2}(q),...\in qk[[q]]$
\end{theorem}

 Conversely we will prove:

 \begin{theorem}
\label{infinitesums}
Any element of the form
$$\overline{f}=\varphi_0(q)+\sum_{s\geq 0} (V^{s+1}(\varphi_{p^s}(q)))\left( \frac{q'}{q^p}\right)^{p^s}\in k((q))[q']$$
with
$\varphi_0(q),\varphi_1(q),\varphi_p(q),\varphi_{p^2}(q),...\in k((q))$
is Laurent $\d-p$-symmetric and
$$``pU" \overline{f}=-\sum_{s\geq 0} V^s(\varphi^{(-1)}_{p^s}(q))+\sum_{s\geq 0} (V^{s+1}(U( \varphi_{p^s}(q))))\left(\frac{q'}{q^p}\right)^{p^s}.$$
If in addition  $\overline{f}\in k[[q]][q']$
(i.e. if
$
\varphi_0(q) \in  k[[q]]$ and $\varphi_1(q),\varphi_p(q),\varphi_{p^2}(q),...\in qk[[q]])$)
 then $\overline{f}$ is Taylor $\d-p$-symmetric.
\end{theorem}

\medskip

\begin{corollary}
\label{voce2}
Let $\overline{f}\in k((q))[q']$ be Laurent $\d-p$-symmetric and let $\overline{\lambda}_p\in k$. Then  $``pT_{\k}(p)"\overline{f}=\overline{\lambda}_p\cdot \overline{f}$ if and only if:

1) $U(\varphi_{p^s}(q))=\overline{\lambda}_p \cdot \varphi_{p^s}(q)$ for all $s\geq 0$ and

2) $\ee\cdot V(\varphi_0(q))-\sum_{s\geq 0} V^s(\varphi_{p^s}^{(-1)}(q))=\overline{\lambda}_p\cdot \varphi_0(q)$.
\end{corollary}

\begin{corollary}
If $\overline{f}\in k[[q]][q']$ is Taylor $\d-p$-symmetric then the series
$``pU"\overline{f}$ and $``pT_{\k}(p)"\overline{f}$ are again Taylor $\d-p$-symmetric.
\end{corollary}

\begin{remark}
It is tempting to conjecture that any Taylor $\d-p$-symmetric series in $k[[q]][q',...,q^{(r)}]$
must belong to $k[[q]][q']$.
\end{remark}

We will first prove  Theorem \ref{infinitesums}.
The plan will be to  first prove this Theorem in case $\overline{f}$ is a monomial
in $k[q,q']$;
cf. Lemma \ref{echelon2} below. This will imply, of course, that Theorem \ref{infinitesums} holds in case
$\overline{f}$ is a finite sum of monomials. The rest of the proof will be devoted to extending
the result from finite to infinite sums of monomials; this will require an analysis
of $(s_1,...,s_p)$-adic convergence of certain series.

\begin{lemma}
\label{echelon2}
For any $n \in {\mathbb Z}$ and $s \in {\mathbb Z}_{\geq 0}$  the element
$$\overline{f}=q^{np^{s+1}}(q')^{p^s}=q^{(n+1)p^{s+1}} \left(\frac{q'}{q^p}\right)^{p^s}\in k((q))[q']$$
  is  Laurent $\d-p$-symmetric (and actually
  Taylor $\d-p$-symmetric if $n \geq 0$.) Moreover
$$
``pU"\overline{f}=\begin{cases}
q^{(n+1)p^s}\left(\frac{q'}{q^p}\right)^{p^s}\ \ \text{if}\ \ p|n+1\\
\ \\
-\frac{q^{(n+1)p^s}}{n+1}\ \ \text{if}\ \ p\not| n+1
\end{cases}
$$
\end{lemma}

{\it Proof}. It is enough to consider the case $s=0$; the general case follows by applying the $p$-power Frobenius.

For $n=-1$ note that
$$q^{-p}q' \equiv \Psi\ \ mod\ \ (p)$$
and so $q^{-p}q'$ is Laurent $\d-p$-symmetric  because $\Psi$ is Laurent $\d-p$-symmetric.
Also $``pU" \overline{f}=\overline{f}$ because $pU\Psi=\Psi$.

Assume now $n\neq -1$.
We have
$$\begin{array}{rcl}
\d(q^{n+1}) & = & \frac{1}{p}[(q^p+pq')^{n+1}-q^{p(n+1)}]\\
\  & \  & \ \\
\  & = & \frac{1}{p}\left[p(n+1) q^{pn} q' +\sum_{j\geq 2} \frac{p^j}{j!} (n+1)...(n-j+2) q^{p(n+1-j)}(q')^j\right]\end{array}$$
For $j \geq 2$ (and since $p \geq 5$) we have
$$v_p\left( \frac{p^j}{j!} \right) \geq j-v_p(j!)\geq j-\frac{j}{p-1}>1.$$
It follows that
\begin{equation}
\label{deltadeputere}
\d(q^{n+1})=(n+1)[q^{pn}q'+pF_{n+1}(q,q')],\ \ \
F_{n+1}(q,q') \in
R[q,q^{-1},q'].
\end{equation}
In particular $\d(q^{n+1})$ is divisible by $n+1$ in $R((q))\h[q']\h$ and we have the following congruence in $R((q))\h[q']\h$:
\begin{equation}
\label{adele}
\frac{1}{n+1}\d(q^{n+1})\equiv q^{np}q'\ \ mod\ \ (p).\end{equation}
By Remark \ref{multe}, assertions 4) and 5), the left hand side of the latter congruence is Laurent $\d-p$-symmetric (and also Taylor
 $\d-p$-symmetric if $n \geq 0$) and hence $q^{pn}q'$ is
 Laurent $\d-p$-symmetric (and also Taylor
 $\d-p$-symmetric if $n \geq 0$).

 To compute $``pU"\overline{f}$ start with the following computation in $R((q))\h[q']\h$:
 $$\begin{array}{rcl}
 p^2(n+1)U\left(\frac{\d(q^{n+1})}{n+1}\right) & = & pU(p\d(q^{n+1})) \\
 \  & \  & \  \\
 \  & = & pU(\phi(q^{n+1}))-pU(q^{p(n+1)})\\
 \  & \  & \  \\
 \  & = & \phi(pU(q^{n+1}))-pU(q^{p(n+1)})\\
\  & \  & \  \\
\  & = & \begin{cases}
-pq^{n+1}\ \ \text{if}\ \ p\not| n+1\\
p\phi(q^{\frac{n+1}{p}})- pq^{n+1}\ \ \text{if}\ \ p|n+1\end{cases}\\
\  & \  & \  \\
\  & = &
\begin{cases}
-pq^{n+1}\ \ \text{if}\ \ p\not| n+1\\
p^2\d(q^{\frac{n+1}{p}})\ \ \text{if}\ \ p|n+1\end{cases}\\
\  & \  & \  \\
\  & = & \begin{cases}
-pq^{n+1}\ \ \text{if}\ \ p\not| n+1\\
p^2 \frac{n+1}{p} \left[ q^{p(\frac{n+1}{p}-1)}q'+pF_{\frac{n+1}{p}}(q,q')\right]\ \ \text{if}\ \ p|n+1\end{cases}
\end{array}$$
from which we get the following congruences mod $p$ in $R((q))\h[q']\h$:
$$ pU(q^{pn}q')\equiv pU \left( \frac{\d(q^{n+1})}{n+1} \right)
\equiv \begin{cases}
-\frac{q^{n+1}}{n+1}\ \ \text{if}\ \ p\not|n+1\\
\ \\
q^{n+1-p}q'\ \ \text{if}\ \ p|n+1.\end{cases}
$$
and we are done.
\qed
\bigskip

\begin{lemma}
\label{ast}
Consider the polynomials
$$s_1,...,s_p,s'_1,...,s'_p, D \in k[q_1,...,q_p,q'_1,...,q_p],\ \ D:=\prod_{i<j}(q_i-q_j).$$
Then the polynomials
$$D^p q'_1,...,D^p q'_p$$
are linear combinations of
$$1,s'_1,...,s'_p$$
with coefficients in $k[q_1,...,q_p]$.
\end{lemma}

{\it Proof}.
For $j=1,...,p$ let $s_{ij}$ be obtained from $s_i$ by setting $q_j=0$; so $s_{ij}$ is the $i$th fundamental symmetric polynomial
in $\{q_1,...,q_p\}\backslash \{q_j\}$. Taking $\d$ in the equalities
$$q_1+...+q_p=s_1,...,q_1...q_p=s_p$$
in $R[q_1,...,q_p,q'_1,...,q'_p]$ and reducing mod $p$
we get the following equalities in $k[q_1,...,q_p,q'_1,...,q'_p]$:
$$
\begin{array}{rcl}
q'_1+...+q'_p & = & s'_1-\gamma_1\\
\  & \  & \  \\
s_{11}^p q'_1+...+s_{1p}^pq'_p & = & s'_2-\gamma_2\\
\  & \  & \  \\
............................... & \  & \  \\
\  & \  & \  \\
s_{p-1,1}^pq'_1+...+s_{p-1,p}^pq'_p & = & s_p'-\gamma_p
\end{array}$$
for some $\gamma_1,...,\gamma_p \in k[q_1,...,q_p]$. View this as a linear system of equations with unknowns
$q'_1,...,q'_p$. We shall be done if we prove that the determinant of the matrix of this system is $\pm D^p$. This follows by taking determinants
in the obvious identity of matrices

\bigskip

$$\left( \begin{array}{rrrr}
q_1^{p-1} & -q_1^{p-2} & ... & 1\\
\  & \  & \  & \  \\
q_2^{p-1} & -q_2^{p-2} & ... & 1\\
\  & \  & \  & \  \\
... & \  & \  & \  \\
\  & \  & \  & \  \\
q_p^{p-1} & -q_p^{p-2} & ... & 1 \end{array}\right)
\left( \begin{array}{rrrr}
1 & 1 & ... & 1\\
\  & \  & \  & \  \\
s_{11} & s_{12} & ... & s_{1p}\\
\  & \  & \  & \  \\
... & \  & \  & \  \\
\  & \  & \  & \  \\
s_{p-1,1} & s_{p-1,2} & ... & s_{p-1,p}\end{array}\right)=(D_{ij})$$

\bigskip

\noindent where
$$D_{ij}=\prod_{s \neq j}(q_i-q_s)$$
and noting that $(D_{ij})$ is a diagonal matrix with determinant $D^2$.
\qed

\begin{lemma}
\label{asta}
Assume the notation of Lemma \ref{ast} and $n \geq 0$. Then the element
$$\sum_{i=1}^p q_i^{np}q'_i \in k[[q_1,...,q_p]][q'_1,...,q'_p]$$
is a  linear combination of
$$1,s'_1,...,s'_p$$
with coefficients in the ideal
$$(s_1,...,s_p)^{[(n+1)/p]-1}k[s_1,...,s_p].$$
\end{lemma}

{\it Proof}.
By Lemma \ref{ast} we can write
$$\sum_{i=1}^p q_i^{np}q'_i=A_0+\sum_{j=1}^p A_j s'_j$$
where $A_j \in k[q_1,...,q_p,D^{-1}]$ for $j=0,...,p$. On the other hand, by (\ref{adele})
$\sum_{i=1}^p q_i^{np}q'_i$ is the reduction mod $p$ of
$$\frac{1}{n+1}\sum_{i=1}^p \d(q_i^{n+1}) \in R[q_1,...,q_p,q'_1,...,q_p].$$
We claim that the following holds:
\begin{equation}
\label{claim}
\sum_{i=1}^p \d(q_i^{n+1}) \in (s_1,...,s_p,s'_1,...,s'_p)^{[(n+1)/p]}R[s_1,...,s_p,s'_1,...,s'_p].
\end{equation}
Assuming (\ref{claim}) is true let us show how to conclude the proof of the Lemma. By (\ref{claim}) we get that
$$\sum_{i=1}^p q_i^{np}q'_i \in (s_1,...,s_p,s'_1,...,s_p)^{[(n+1)/p]}k[s_1,...,s_p,s'_1,...,s'_p].$$
So we have
$$\sum_{i=1}^p q_i^{np}q'_i=\sum B_{i_1...i_p}(s'_1)^{i_1}...(s'_p)^{i_p}$$
where
$$B_{i_1...i_p} \in (s_1,...,s_p)^{[(n+1)/p]-i_1-...-i_p}k[s_1,...,s_p].$$
Since $s'_1,...,s'_p$ are algebraically independent over $k[q_1,...,q_p]$ we get
$$\begin{array}{rcl}
A_0 & = & B_{0...0}\\
\  & \  & \  \\
A_1 & = & B_{10...0}\\
\  & \  & \  \\
A_2 & = & B_{010...0},\ \text{etc}\end{array}$$
hence
$$A_j \in (s_1,...,s_p)^{[(n+1)/p]-1}k[s_1,...,s_p],\ \ j=0,...,p$$
which ends the proof of the Lemma.

To check (\ref{claim}) above note that
$$\sum_{i=1}^p \d(q_i^{n+1})=\d\left(\sum_{i=1}^p q_i^{n+1}\right)+\frac{\left(\sum_{i=1}^p q_i^{n+1}\right)^p-\sum_{i=1}^p q_i^{(n+1)p}}{p}.$$
The second term in the right hand side of the above equation is a homogeneous polynomial in $q_1,...,q_p$ of degree $(n+1)p$ hence it is a weighted homogeneous polynomial in $s_1,...,s_p$ of weight $(n+1)p$
where $s_1,...,s_p$ are given weights $1,...,p$ respectively. Hence this polynomial is a sum of monomials
in $s_1,...,s_p$ of degree $\geq n+1$. Similarly $\sum_{i=1}^p q_i^{n+1}$ is a sum of monomials in $s_1,...,s_p$ of degree $\geq [(n+1)/p]$. This implies that $\d(\sum_{i=1}^p q_i^{n+1})$
is a sum of monomials in $s_1,...,s_p,s'_1,...,s'_p$ of degree $\geq [(n+1)/p]$ which proves (\ref{claim}).
\qed

\medskip

{\it Proof of Theorem \ref{infinitesums}.}
In view of Lemma \ref{echelon2} (which treats the case of monomials) we see that
in order to prove that  $\overline{f}$ in the statement of the Theorem is Laurent (respectively Taylor) $\d-p$-symmetric it is enough to show that any series of the form
$$\overline{f}=\sum_{n=0}^{\infty} \overline{c}_n q^{pn}q' \in k[[q]][q']$$
is Taylor $\d-p$-symmetric. By Lemma \ref{asta} we may write
$$\sum_{i=1}^p q_i^{np}q'_i=G_{0n}+\sum_{j=1}^p G_{jn} s'_j$$
where
$$G_{jn}\in (s_1,...,s_p)^{[(n+1)/p]-1}k[s_1,...,s_p],\ \ j=0,...,p.$$
Since $G_j:=\sum_{n=0}^{\infty} \overline{c}_n G_{jn}$ are convergent in $k[[s_1,...,s_p]]$ we have
$$\sum_{i=1}^p \overline{f}(q_i)=G_0+\sum_{j=1}^p G_j s'_j \in k[[s_1,...,s_p]][s'_1,...,s'_p]$$
which proves that $\overline{f}$ is Taylor $\d-p$-symmetric. The assertion about $``pU" \overline{f}$ follows from Lemma \ref{echelon} by taking limits.
\qed

Next we proceed to proving Theorem \ref{conv}.
We need a preparation.
Let $C_p(q_1,q_2):=\frac{q_1^p+q_2^p-(q_1+q_2)^p}{p}\in {\mathbb Z}[q_1,q_2]$.
We start with a version of Lemma
\ref{ast}:

\begin{lemma}
\label{shape}
Consider the elements  $\sigma=q_1+q_2\in k[q_1,q_2]$ and $\pi=q_1q_2\in k[q_1,q_2]$ and let $\gamma \in k[q_1,q_2]$
be the image of $C_p(q_1,q_2)\in {\mathbb Z}[q_1,q_2]$. Then
$$q'_1=\frac{\pi'-q_1^p\sigma'+q_1^p\gamma}{(q_2-q_1)^p},\ \ q'_2=-\frac{\pi'-q_2^p\sigma'+q_2^p\gamma}{(q_2-q_1)^p}$$
in the ring
$$k[q_1,q_2,q'_1,q'_2,\frac{1}{q_2-q_1}].$$
\end{lemma}

{\it Proof}.
Applying $\d$ to the defining equations of $\sigma$ and $\pi$ we get
$$\begin{array}{rcl}
q_1'+q_2' & = & \sigma'-\gamma\\
\  & \  & \  \\
q_2^pq'_1+ q_1^p q'_2 & = & \pi'\\
\end{array}$$
and solve for $q'_1,q_2'$.
\qed

\bigskip

For the next Lemma let us denote by $v_{q_2-q_1}:k((q_1,q_2))^{\times}\ra {\mathbb Z}$ the normalized
valuation on the fraction field $k((q_1,q_2))$ of $k[[q_1,q_2]]$ attached to the irreducible series $q_2-q_1 \in k[[q_1,q_2]]$; in other words,
if $0\neq F(q_1,q_2)\in k[[q_1,q_2]]$ then $v_{q_2-q_1}(F)$ is the maximum integer $i$ such that
$(q_2-q_1)^i$ divides $F$ in $k[[q_1,q_2]]$.

\begin{lemma}
\label{val}
Let $\Phi(q)=\sum_{m=0}^{\infty} \beta_mq^m \in k[[q]]$, $\Phi \not\in k$, $Supp\ \Phi:=\{m \in {\mathbb Z}_{\geq 0};\beta_m\neq 0\}$.
Then
$$v_{q_2-q_1}(\Phi(q_2)-\Phi(q_1))=p^{\min\{v_p(m);0\neq m \in Supp\ \Phi\}}.$$
\end{lemma}

{\it Proof}.
We have
$$\begin{array}{rcl}
\Phi(q_2)-\Phi(q_1) & = & \sum_{(n,p)=1}\sum_{i=0}^{\infty} \beta_{np^i}(q_2^{np^i}-q_1^{np^i})\\
\  & \  & \  \\
\  & = &  \sum_{i=0}^{\infty} (q_2-q_1)^{p^i}G(q_1,q_2)\end{array}$$
where
$$G_i(q_1,q_2)=\sum_{(n,p)=1} \beta_{np^i}(q_2^{(n-1)p^i}+q_2^{(n-2)p^i}q_1^{p^i}+...+q_1^{(n-1)p^i}).$$
Let $i_0=\min\{v_p(m);0\neq m \in Supp\ \Phi\}$. Then $\beta_{np^{i}}=0$ for all $(n,p)=1$ and $i<i_0$ and there exists $n_0$, $(n_0,p)=1$
such that $\beta_{n_0p^{i_0}}\neq 0$. It is enough to show that
$G_{i_0}(q_1,q_2)$ is not divisible by $q_2-q_1$ in $k[[q_1,q_2]]$ equivalently that $G(q,q)\neq 0$. But
$$G_{i_0}(q,q)=\sum_{(n,p)=1} n \beta_{np^{i_0}} q^{(n-1)p^{i_0}}\neq 0.$$
\qed

\bigskip

{\it Proof of Theorem \ref{conv}}.
We proceed by induction on the degree $\deg(\overline{f})$ of $\overline{f}$ viewed as a polynomial in $q'$ with coefficients in $k[[q]]$.
If this degree is $0$ we are done. Assume now the degree is $\geq 1$. We may assume $\overline{f}(0,0)=0$.

By hypothesis,
$$\overline{f}(q_1,q'_1)+...+\overline{f}(q_p,q'_p)=G$$
in $k[[q_1,...,q_p]][q'_1,...,q'_p]$, where $G\in k[[s_1,...,s_p]][s'_1,...,s'_p]$.
Setting $q_3=...=q_p=0$ and $q'_3=...=q'_p=0$ we get
\begin{equation}
\label{cucu}
\overline{f}(q_1,q'_1)+\overline{f}(q_2,q'_2)=G(\sigma,\pi,0,...,0,\sigma',\pi',0,...,0).\end{equation}
Note that $k[[q_1,q_2]]$ is a finite $k[[\sigma,\pi]]$-algebra so $\sigma',\pi'$ are algebraically independent over $k((q_1,q_2))$.
By Lemma \ref{shape} the left hand side of (\ref{cucu}) is a  polynomial $H$ in $\sigma', \pi'$ with coefficients in $k((q_1,q_2))$. On the other hand since $H$ is in the right hand side of (\ref{cucu})
$H$ has  coefficients in $k[[q_1,q_2]]$. Hence each non-zero coefficient  of the polynomial $H$ has $v_{q_2-q_1}$-adic valuation $\geq 0$.
 Now write
$$\overline{f}(q,q')=\sum_{m'}\Phi_{m'}(q)(q')^{m'},\ \ \ \Phi_{m'} \in k[[q]].$$
Also write each $m'$ as
 $m'=n'p^{i'}$ with $n'$ not divisible by $p$. Using Lemma \ref{shape} we have $H=\sum_{m'} H_{m'}$ where
\begin{equation}
\label{FF}
H_{m'}=\frac{F_{m'}}{(q_2-q_1)^{n'p^{i'+1}}}\end{equation}
where $F_{m'} \in k((q_1,q_2))[\sigma',\pi']$ is given by
$$\begin{array}{rcl}
F_{m'} & = &
\Phi_{m'}(q_1)\left((\pi')^{p^{i'}}-q_1^{p^{i'+1}}(\sigma')^{p^{i'}}+q_1^{p^{i'+1}}\gamma^{p^{i'}}\right)^{n'}\\
\  & \  & \  \\
\  & \  &
+(-1)^{n'}\Phi_{m'}(q_2)
\left((\pi')^{p^{i'}}-q_2^{p^{i'+1}}(\sigma')^{p^{i'}}+q_2^{p^{i'+1}}\gamma^{p^{i'}}\right)^{n'}.\end{array}
$$
Note that the coefficient of $(\pi')^{m'}$ in $F_{m'}$ is
\begin{equation}
\label{star}
\Phi_{m'}(q_1)+(-1)^{n'}\Phi_{m'}(q_2)
\end{equation}
 while the coefficient of $(\pi')^{m'-p^{i'}}(\sigma')^{p^{i'}}$ in $F_{m'}$ is
 \begin{equation}
 \label{starr}
 -n'\left( q_1^{p^{i'+1}}\Phi_{m'}(q_1)+(-1)^{n'}q_2^{p^{i'+1}}\Phi_{m'}(q_2)\right).
 \end{equation}
Let now $m'=\deg(\overline{f})$.
 If $n'$ is even the polynomial (\ref{star})  has $v_{q_2-q_1}$-adic valuation $0$ which contradicts
 the fact that the non-zero coefficients of $H$ have $v_{q_2-q_1}$-adic valuation $\geq 0$.
So $n'$ is odd. By Lemma \ref{val}
the $v_{q_2-q_1}$-adic valuation of (\ref{star}) equals
$$p^{\min\{v_p(m);0 \neq m \in Supp\ \Phi_{m'}\}},\ \ \ \text{if}\ \ \Phi_{m'}\not\in k.$$
Also the
$v_{q_2-q_1}$-adic valuation of (\ref{starr}) equals
$$p^{\min\{v_p(m);m\in Supp(q^{p^{i'+1}}\Phi_{m'})\}}=p^{\min\{v_p(m+p^{i'+1});m \in Supp\ \Phi_{m'}\}}.$$
 By the fact that the non-zero coefficients of $H$ have $v_{q_2-q_1}$-adic valuation $\geq 0$ we get that
\begin{equation}
\label{apple}
p^{\min\{v_p(m);0 \neq m \in Supp\ \Phi_{m'}\}} \geq n'p^{i'+1}\ \ \text{if}\ \ \Phi_{m'}\not\in k
\end{equation}
and
\begin{equation}
\label{applee}
  p^{\min\{v_p(m+p^{i'+1});m \in Supp\ \Phi_{m'}\}}  \geq n'p^{i'+1}.
\end{equation}
From (\ref{apple}) we get
\begin{equation}
\label{conc1}
v_p(m) \geq i'+1\ \ \text{for all} \ \ 0\neq m \in Supp\ \Phi_{m'},\ \  \text{if}\ \ \Phi_{m'}\not\in k.
\end{equation}
We claim now  that $n'=1$. Assume $n'\geq 2$. By (\ref{apple})
$$v_p(m)>i'+1\ \ \text{for all}\ \ 0 \neq m \in Supp\ \Phi_{m'},\ \ \text{if}\ \ \Phi_{m'}\not\in k.$$
Hence
$$v_p(m+p^{i'+1})=i'+1\ \ \text{for all}\ \ m \in Supp\ \Phi_{m'}.$$
 By (\ref{applee}) $p^{i'+1}\geq 2p^{i'+1}$, a contradiction.
This ends the proof that $n'=1$.

By (\ref{conc1})
$$\Phi_{m'}(q)(q')^{m'}=(V^{i'+1}\varphi)(q')^{p^{i'}}$$
for some $\varphi \in k[[q]]$. By Lemma \ref{echelon2} $\Phi_{m'}(q)(q')^{m'}$
is Taylor $\d-p$-symmetric hence so is $\overline{f}-\Phi_{m'}(q)(q')^{m'}$ which has smaller
degree than $\overline{f}$. We conclude by the induction hypothesis.
\qed

\section{Multiplicity one}

We begin by recalling the well known situation for series in $k[[q]]$.
Then we proceed with our main results about $\d$-series in $k[[q]][q']$.

Throughout this section we fix $\k \in {\mathbb Z}_{\geq 0}$.

\begin{definition}
A series $\varphi\in qk[[q]]$ is said to be an {\it eigenvector of all Hecke operators
$T_{\k+2}(n)$, $T_{\k+2}(p)$, $(n,p)=1$}, with eigenvalues $\overline{\lambda}_n, \overline{\lambda}_p \in k$
if
 $\varphi\neq 0$ and the following hold:
\begin{equation}
\label{eqvarphi}
\begin{cases}
T_{\k+2}(n)\varphi  =  \overline{\lambda}_n \cdot \varphi, \ \ \  (n,p)=1\\
\  \\
T_{\k+2}(p)\varphi  =  \overline{\lambda}_p\cdot \varphi.\end{cases}
\end{equation}
Of course the last equation in (\ref{eqvarphi}) is equivalent to
$$U\varphi=\overline{\lambda}_p\cdot \varphi.$$
\end{definition}

\begin{proposition}
Assume $\varphi\in qk[[q]]$ is an eigenvector of all Hecke operators
$T_{\k+2}(n)$, $T_{\k+2}(p)$, $(n,p)=1$, with eigenvalues $\overline{\lambda}_n, \overline{\lambda}_p \in k$
Then there exists $\gamma \in k^{\times}$ such that
\begin{equation}
\label{defoffi}
\varphi(q):=\gamma \cdot
\sum_{(n,p)=1} \sum_{i \geq 0} \overline{\lambda}_n\overline{\lambda}_p^i \cdot q^{np^i}.\end{equation}
\end{proposition}

{\it Proof}.
Pick out coefficient of $q$ in the first equation (\ref{eqvarphi}) and the coefficient of $q^m$, $m \geq 1$ in the second equation (\ref{eqvarphi}). (Here we use the convention that $0^0=1$.)
\qed

\bigskip

 \begin{definition}
  A $\d$-series  $\overline{f}=\overline{f}(q,q')\in k[[q]][q']$
  is said to be an {\it eigenvector of all Hecke operators
$nT_{\k}(n)$, ``$pT_{\k}(p)$", $(n,p)=1$}, with eigenvalues $\overline{\lambda}_n, \overline{\lambda}_p \in k$
if $\overline{f}$
  is Taylor $\d-p$-symmetric and satisfies
 \begin{equation}
 \label{eqf}
 \begin{cases}
 nT_{\k}(n)\overline{f}  =  \overline{\lambda}_n \cdot \overline{f}, \ \ \  (n,p)=1;\\
 \  \\
``pT_{\k}(p)"\overline{f}  =  \overline{\lambda}_p \cdot \overline{f}. \end{cases}
\end{equation}
 \end{definition}

\begin{theorem}
\label{arnab}
 Assume  $\overline{f}=\overline{f}(q,q')\in k[[q]][q']$,  $\overline{f}\not\in k$, is
 an  eigenvector of all Hecke operators
$nT_{\k}(n)$, ``$pT_{\k}(p)$", $(n,p)=1$, with eigenvalues $\overline{\lambda}_n, \overline{\lambda}_p \in k$.
Then there exists  $\varphi=\varphi(q)\in qk[[q]]$ and  $c, c_i \in k$, $i \geq 0$,
with $\ee \cdot c_{i-1}=\overline{\lambda}_p\cdot c_i$ for $i\gg 0$,
such that
$\varphi$ is an eigenvector of all Hecke operators
$T_{\k+2}(n)$, $T_{\k+2}(p)$, $(n,p)=1$, with the same eigenvalues $\overline{\lambda}_n, \overline{\lambda}_p$ and such that
\begin{equation}
\label{zor}
\begin{array}{rcl}
\overline{f} & = & c+
\left(\sum_{i \geq 0} c_i F_{/k}^i\right)\varphi^{\sharp,2},\\
\  & \  & \  \\
\varphi^{\sharp,2} & := & \varphi^{(-1)}-\overline{\lambda}_p \cdot V(\varphi)
\frac{q'}{q^p}+\ee\cdot V^2(\varphi)\left(\frac{q'}{q^p}\right)^p.\end{array}
\end{equation}
\end{theorem}

\begin{remark}
One can also write $\overline{f}$ in (\ref{zor}) as
$$\begin{array}{rcl}
\overline{f} & = & c+
\sum_{i \geq 0} c_i \left[ V^i(\varphi^{(-1)})-\overline{\lambda}_p \cdot V^{i+1}(\varphi)
\left(\frac{q'}{q^p}\right)^{p^i}+\ee\cdot V^{i+2}(\varphi)\left(\frac{q'}{q^p}\right)^{p^{i+1}}\right]\\
\  & \  & \  \\
\  & = & c+\left( \sum_{i\geq 0} c_iV^i\right) \varphi^{(-1)} + \sum_{i \geq 0}
(\ee c_{i-1}-\overline{\lambda}_p c_i) V^{i+1}(\varphi) \left(\frac{q'}{q^p}\right)^{p^i},\end{array}$$
where $c_{-1}:=0$.
Note that the condition that $\ee \cdot c_{i-1}=\overline{\lambda}_p\cdot c_i$ for $i\gg 0$
insures  that the right hand side of the first equation in (\ref{zor}) is a polynomial in the variable $q'$.
\end{remark}

\begin{remark} Looking at the constant terms in (\ref{eqf}) one sees  that if $c\neq 0$
then
\begin{equation}
\label{oprit}
\begin{cases}
\overline{\lambda}_n  =  n \cdot \sum_{A|n} \epsilon(A)A^{\k-1}, \ \  (n,p)=1;\\
\ \\
\overline{\lambda}_p  =  \ee.  \ \end{cases}\end{equation}
\end{remark}

Conversely we will prove:

\begin{theorem}
\label{curcan}
Let $\k \in {\mathbb Z}_{\geq 0}$.
 Assume $\varphi=\varphi(q)\in qk[[q]]$
is an eigenvector of all Hecke operators $T_{\k+2}(n)$, $T_{\k+2}(p)$, $(n,p)=1$, with eigenvalues
  $\overline{\lambda}_n,\overline{\lambda}_p \in k$.
  Let  $c_i \in k$ for $i \geq 0$ with $\ee \cdot c_{i-1}=\overline{\lambda}_p\cdot c_i$ for $i\gg 0$.
  Also let $c$ be an arbitrary element in $k$ or $0$ according as equations (\ref{oprit}) hold or fail respectively.
  Let $\overline{f}\in k[[q]][q']$ be defined by
  Equation (\ref{zor}).
 Then $\overline{f}$ an  eigenvector of all Hecke operators
$nT_{\k}(n)$, ``$pT_{\k}(p)$", $(n,p)=1$, with the same eigenvalues $\overline{\lambda}_n, \overline{\lambda}_p$.
\end{theorem}

Let $k[F_{/k}]$ be the $k$-algebra generated by $F_{/k}$ which is a commutative polynomial ring in one variable. Note that $k[[q]][q']$ is a $k[F_{/k}]$-module and the $k$-linear space of series $\overline{f}(q,q')\in k[[q]][q']$ with $f(0,0)=0$
is a torsion free $k[F_{/k}]$-submodule.
Note also that the ideal $qk[[q]]$ is a torsion free module over the ring $k[[F_{/k}]]$ of power series in $F_{/k}$.
Finally recall that a $\d$-series $\overline{f}(q,q')\in k[[q]][q']$ is called {\it primitive}
if $U(\overline{f}(q,0))=0$.
 Theorems \ref{arnab} and \ref{curcan} immediately imply:

\begin{corollary}
 Fix   $\overline{\lambda}_n \in k$ for $(n,p)=1$ and  $\overline{\lambda}_p \in k$.
 Let ${\mathcal F}$
 be the $k$-linear space of all the   $\d$-series
$\overline{f}=\overline{f}(q,q')\in k[[q]][q']$ with $f(0,0)=0$
which are either $0$ or are eigenvectors of all Hecke operators
$nT_{\k}(n)$, ``$pT_{\k}(p)$", $(n,p)=1$, with  eigenvalues $\overline{\lambda}_n, \overline{\lambda}_p \in k$.
We have ${\mathcal F}\neq 0$  if and only if there exists an eigenvector $\varphi\in qk[[q]]$
 of all Hecke operators $T_{\k+2}(n)$, $T_{\k+2}(p)$, $(n,p)=1$, with eigenvalues
 $\overline{\lambda}_n, \overline{\lambda}_p$.
 Assume furthermore that this is the case and let
  $\varphi^{\sharp,2}$ be defined as in  (\ref{zor}).
Then $\varphi^{\sharp,2}$ belongs to ${\mathcal F}$ and is a primitive $\d$-series; also any primitive $\d$-series in ${\mathcal F}$
is a $k$-multiple of $\varphi^{\sharp,2}$. Furthemore the following hold:

1) If $\k>0$, $\overline{\lambda}_p=0$ then ${\mathcal F}$ is a free $k[[F_{/k}]]$-submodule of $k[[q]]$ of rank $1$ with basis $\varphi^{\sharp,2}=\varphi^{(-1)}$.

2) If either $\k >0$, $\overline{\lambda}_p\neq 0$ or $\k=0$, $\overline{\lambda}_p=0$ then ${\mathcal F}$ is a free $k[F_{/k}]$-submodule
of $k[[q]][q']$ of rank one with basis $\varphi^{\sharp,2}$.

3) If $\k=0$, $\overline{\lambda}_p\neq 0$ then ${\mathcal F}$ is a free $k[F_{/k}]$-submodule of $k[[q]][q']$ of rank $1$ with basis
\begin{equation}
\label{varphisharpunu}
\varphi^{\sharp,1}:=\left(\sum_{i\geq 0} (\overline{\lambda}_p)^{-i} F^i_{/k}\right)\varphi^{\sharp,2}.\end{equation}
\end{corollary}

\begin{remark}
Note that
$$\varphi^{\sharp,1}=\left(\sum_{i\geq 0} (\overline{\lambda}_p)^{-i} V^i\right)\varphi^{(-1)}-\overline{\lambda}_p\cdot V(\varphi)\cdot \frac{q'}{q^p}$$
and also that  $\varphi^{\sharp,1}$ is the unique element of $qk[[q]]$ satisfying the equation
$$V(\varphi^{\sharp,1})-\overline{\lambda}_p \varphi^{\sharp,1}+\overline{\lambda}_p \varphi^{\sharp,2}=0.$$
\end{remark}

\bigskip

 {\it Proof of Theorem \ref{arnab}}.
 For any series $\beta(q) \in k[[q]]$ write
$$\beta(q)=\sum_{m\geq 0} a_m(\beta) q^m.$$
 By Theorem \ref{conv} and  Corollaries
\ref{voce} and \ref{voce2} $\overline{f}$ has the form
(\ref{blabla}) and

\bigskip

\begin{equation}
\label{urrs}\begin{array}{rcll}
T_{\k}(n)\varphi_0 & = & \frac{\overline{\lambda}_n}{n} \cdot \varphi_0, & (n,p)=1\\
\  & \  & \  & \  \\
T_{\k+2p^s}(n)\varphi_{p^s} & = & \overline{\lambda}_n \cdot \varphi_{p^s}, & (n,p)=1,\ s\geq 0\\
\  & \  & \  & \  \\
U(\varphi_{p^s}) & = & \overline{\lambda}_p \cdot \varphi_{p^s}, & s \geq 0\\
\  & \  & \  & \  \\
\ee\cdot V(\varphi_0) -\sum_{s \geq 0} V^s(\varphi_{p^s}^{(-1)}) & = & \overline{\lambda}_p \cdot \varphi_0. & \  \end{array}
\end{equation}

\bigskip

\noindent In particular the following equalities hold:

\bigskip

\begin{equation}
\label{ocup}\begin{array}{rcll}
a_{np^s}(\varphi_0) & = & \frac{\overline{\lambda}_n}{n} \cdot a_{p^s}(\varphi_0), & (n,p)=1, s\geq 0,\\
\  & \  & \  & \  \\
a_{n}(\varphi_{p^s}) & = & \overline{\lambda}_n  \cdot a_1(\varphi_{p^s}), & (n,p)=1,  s\geq 0,\\
\  & \  & \  & \  \\
a_{mp}(\varphi_{p^s}) & = & \overline{\lambda}_p \cdot a_m(\varphi_{p^s}), & m\geq 1, s\geq 0,\\
\  & \  & \  & \  \\
\ee\cdot a_{p^{s-1}}(\varphi_0)-  a_1(\varphi_{p^s}) & = & \overline{\lambda}_p \cdot a_{p^s}(\varphi_0), & s \geq 0,\end{array}
\end{equation}

\bigskip

\noindent where by convention we set $a_{p^{s-1}}(\varphi_0)=0$ if $s=0$.
Let $c=a_0(\varphi_0)$ and $c_i=a_{p^i}(\varphi_0)$ for $i \geq 0$.
By (\ref{ocup}) we  get
\begin{equation}
\label{fom}
\begin{array}{rcll}
 a_{np^i}(\varphi_0) & = & \frac{\overline{\lambda}_n}{n} \cdot c_i, & (n,p)=1,\ \ i \geq 0\\
 \  & \ & \  & \ \\
 a_{np^i}(\varphi_{p^s}) & = &  \overline{\lambda}_n \overline{\lambda}_p^i \cdot (\ee \cdot c_{s-1}-\overline{\lambda}_p c_s),
& (n,p)=1,\ \ i \geq 0,\ s\geq 0,\end{array}\end{equation}
where $c_{-1}:=0$.
Define $\varphi$ by the equality (\ref{defoffi}) with $\gamma=1$.

Assume first that  there is an $s \geq 0$ such that $a_1(\varphi_{p^s}) \neq 0$. Then
 $\varphi_{p^s}$ is a non-zero multiple of $\varphi$ so
  (\ref{eqvarphi})
follows from (\ref{urrs}) and (\ref{zor}) follows from (\ref{fom}).
Since $\overline{f}$ is a polynomial in $q'$ we get that $\ee \cdot c_{s-1}-\overline{\lambda}_p c_s=0$ for $s\gg 0$.

 Assume now that $a_1(\varphi_{p^s})=0$ for all $s \geq 0$.
Then $\varphi_{p^s}=0$ for all $s \geq 0$ hence $\overline{f}=\varphi_0$. By the last equation in (\ref{urrs})
and since $\varphi_0\not\in k$ we get $\ee=\overline{\lambda}_p=0$.
Then the right hand side of (\ref{zor}) becomes
\begin{equation}
\label{hihi}
c+\sum_{i \geq 0}  \sum_{(n,p)=1} c_i \frac{\overline{\lambda}_n}{n} q^{np^i}.\end{equation}
By the first equation in (\ref{fom}) we get that (\ref{hihi}) equals $\varphi_0=\overline{f}$;
so equation (\ref{zor}) holds. Clearly $U\varphi=0$ so the second equality in (\ref{eqvarphi}) holds. Finally, since $\varphi_0 \not\in k$ we may
 write $\varphi_0=F^d_{/k}\tilde{\varphi}_0$ with $\tilde{\varphi}_0\in k[[q]]$ and $d$ maximal with this property; in particular $c_d \neq 0$.
Note that $\theta \tilde{\varphi}_0=c_d \varphi$. Also
by (\ref{urrs}) we have $T_{\k}(n)\tilde{\varphi}_0=\frac{\overline{\lambda}_n}{n} \tilde{\varphi}_0$ for $(n,p)=1$. Hence
$$T_{\k+2}(n)\varphi=c_d^{-1} T_{\k+2}(n)\theta \tilde{\varphi}_0=c_d^{-1} n \theta (T_{\k}(n)\tilde{\varphi}_0)=c_d^{-1} \overline{\lambda}_n \theta \tilde{\varphi}_0=\overline{\lambda}_n \varphi$$
and so the first equality in (\ref{eqvarphi}) holds. This ends the proof.
\qed

\bigskip

{\it Proof of Theorem \ref{curcan}}.
This follows directly from Corollary \ref{buium} and Theorem \ref{infinitesums} using the following facts (which are direct consequences of the formulae for  the Hecke operators acting on Fourier coefficients  (\ref{cuani})):
$$\begin{array}{rcll}
T_{\k+2p^i}(n)\varphi & = & \overline{\lambda}_n \cdot \varphi, & (n,p)=1, \ \ \ i\geq 0\\
\  & \  & \ & \  \\
T_{\k}(n)(\varphi^{(-1)}) & = & \frac{\overline{\lambda}_n}{n} \cdot \varphi, & (n,p)=1.  \end{array}$$
\qed

\section{$\d$-modular forms}

\subsection{Review of classical modular forms}
Start by recalling some basic facts about modular forms; cf.  \cite{DI}.
Let $N>4$ be an integer and let $B$ be a ${\mathbb Z}[1/N,\zeta_N]$-algebra.
Let $Y=Y_1(N)$ be the affine modular curve over $B$ classifying pairs $(E,\alpha)$ consisting of elliptic curves $E$ over $B$-algebras plus a
 level $\Gamma_1(N)$ structure
$\alpha:{\mathbb Z}/N{\mathbb Z} \ra E$. Let $Y_{ord}$ be the ordinary locus in $Y$ (i.e. the locus where the Eisenstein form $E_{p-1}$
is invertible).  Let $X$ be $Y$ or $Y_{ord}$. Let $L$ be the line bundle on $X$, direct image of
 the sheaf of relative differentials on the universal elliptic curve over $X$, and let
\begin{equation}
\label{ve}
V=Spec\left(\bigoplus_{\k \in {\mathbb Z}}L^{\otimes \k}\right)\ra X\end{equation}
 be the  ${\mathbb G}_m$-torsor associated to $L$.

Set $M=\cO(V)=\bigoplus_{\k \in {\mathbb Z}}L^{\otimes \k}$. Recall that there is a  Fourier expansion map $$E:M\ra B((q))$$
 defined by the cusp $\Gamma_1(N)\cdot \infty$ \cite{DI}, p. 112. Recall
also that $Y$ has a natural compactification, $X_1(N)$, equipped with a natural line bundle, still denoted by $L$,
  extending the line bundle $L$ on $Y$, such that the space of classical modular forms, $M(\Gamma_1(N),B,\k)\subset L^{\otimes \k}$, on $\Gamma_1(N)$ of weight $\k$, defined over $B$ identifies with $H^0(X_1(N),L^{\otimes \k})$. Recall that the diamond operators act on
  $M(\Gamma_1(N),B,\k)$; the invariant elements form the space $M(\Gamma_0(N),B,\k)$ of classical modular forms on $\Gamma_0(N)$ of weight $\k$ defined over $B$. Recall  the $q$-expansion principle:
  for any $B$ as above there is an induced injective Fourier expansion map $E:M(\Gamma_1(N),B,\k)\ra B[[q]]$ and if $B' \subset B$ then $M(\Gamma_1(N),B',\k)$ identifies with the group of all $f \in M(\Gamma_1(N),B,\k)$ such that $E(f)\in B'[[q]]$. Recall also the following base change property: if $B'$
  is any $B$-algebra and either $B'$  is flat over $B$ or  $\k\geq 2$ and $N$ is invertible in $B'$
then the map $M(\Gamma_1(N),B,\k)\otimes_B B'\ra M(\Gamma_1(N),B',\k)$ is an isomorphism; cf. \cite{DI}, p.111.

\subsection{$\d$-series from classical modular forms}
\begin{theorem}
 \label{nopiine}
 Let $\k \in {\mathbb Z}_{\geq 0}$
and let  $f(q)=\sum_{m\geq 1}a_mq^m\in q{\mathbb Z}_p[[q]]$ be a series satisfying
   $a_1=1$ and
   \begin{equation}
   \label{knacond}
   \begin{cases}
   a_{p^i n}  =  a_{p^i} a_n\ \ \text{for $(n,p)=1,\ i\geq 0$}\\
   a_{p^{i-1}}a_p  =  a_{p^i}+p^{\k+1}a_{p^{i-1}}\ \ \text{for $i \geq 2$}.\end{cases}
   \end{equation}
Let $\varphi:=\overline{f}=\sum_{m \geq 1} \overline{a}_m q^m \in q{\mathbb F}_p[[q]]$ be the reduction mod $p$ of $f(q)$. Then  the series
\begin{equation}
\label{inpat}
f^{\sharp,2}=f^{\sharp,2}(q,q',q''):=\frac{1}{p}
\cdot \sum_{n \geq 1}\frac{a_n}{n}(p^{\k}\phi^2(q)^n-a_p \phi(q)^n+pq^n) \in {\mathbb Q}_p[[q,q',q'']]\end{equation}
belongs to ${\mathbb Z}_p[[q]][q',q'']\h$ and  its reduction mod $p$ equals
\begin{equation}
\label{baie}
\overline{f^{\sharp,2}}=\overline{f^{\sharp,2}(q,q',q'')}=\varphi^{(-1)}-\overline{a}_p V(\varphi)\frac{q'}{q^p}+ \ee \cdot V^{2}(\varphi)\left(\frac{q'}{q^p}\right)^{p}
\in {\mathbb F}_p[[q]][q'].\end{equation}
\end{theorem}

{\it Proof}.
For $\k=0$ the argument is in \cite{BP}; the case $\k>0$ is entirely similar.
(Note that  the form $f_{[a_p]}^{(0)}$ in \cite{BP} is congruent mod $p$ to $f$ itself.)  \qed

\begin{remark}
Note that conditions (\ref{knacond}) imply that $U \varphi=\overline{a}_p \cdot \varphi$.
\end{remark}

\begin{example}
\label{gura}
Let $\k \in {\mathbb Z}_{\geq 0}$ and let $F\subset {\mathbb C}$ be a number field with ring of integers $\cO_F$.
 Let
 \begin{equation}
 \label{cameo}
 f(q)=\sum_{m\geq 1}a_mq^m\in q\cO_F[[q]]\end{equation} be the Fourier expansion of a
 cusp form  $$f\in M(\Gamma_0(N),\cO_F,\k+2).$$
  Assume $a_1=1$ and assume $f(q)$  is an eigenvector for all the Hecke operators
 $T_{\k+2}(n)$ with $n \geq 1$.
 Assume $p$ is a rational prime that splits completely in $F$, consider an embedding $\cO_F \subset {\mathbb Z}_p$, view $f(q)$ as an element of $q{\mathbb Z}_p[[q]]$,
 and let
  $\varphi:=\overline{f}=\sum_{m \geq 1} \overline{a}_m q^m \in q{\mathbb F}_p[[q]]$ is the reduction mod $p$ of $f(q)$.
Then the equalities (\ref{knacond}) hold. So by Theorem  \ref{nopiine} the series
\begin{equation}
\label{inpat}
f^{\sharp,2}=f^{\sharp,2}(q,q',q''):=\frac{1}{p}
\cdot \sum_{n \geq 1}\frac{a_n}{n}(p^{\k}\phi^2(q)^n-a_p \phi(q)^n+pq^n) \in {\mathbb Q}_p[[q,q',q'']]\end{equation}
belongs to ${\mathbb Z}_p[[q]][q',q'']\h$ and  its reduction mod $p$ equals
\begin{equation}
\label{baie}
\overline{f^{\sharp,2}}:=\overline{f^{\sharp,2}(q,q',q'')}=\varphi^{(-1)}-\overline{a}_p V(\varphi)\frac{q'}{q^p}+ \ee \cdot V^{2}(\varphi)\left(\frac{q'}{q^p}\right)^{p}
\in {\mathbb F}_p[[q]][q'].\end{equation}
Note also that $T_{\k+2}(n)\varphi=\overline{a}_n \cdot \varphi$ for $(n,p)=1$ and $U\varphi=\overline{a}_p \cdot \varphi$. So by Theorem \ref{curcan}
   $\overline{f^{\sharp,2}}=\varphi^{\sharp,2}$ is
   an eigenvector of the Hecke operators $nT_{\k}(n)$, ``$pT_{\k}(p)$", $(n,p)=1$, with eigenvalues $\overline{a}_n,\overline{a}_p$.
   Also, by the same Theorem, if in addition $\overline{a}_p\neq 0$ and $\k=0$, then the series $\varphi^{\sharp,1}$ in (\ref{varphisharpunu})
   is also an eigenvector of the Hecke operators $nT_{\k}(n)$, ``$pT_{\k}(p)$", $(n,p)=1$, with eigenvalues $\overline{a}_n,\overline{a}_p$.
\end{example}

\begin{example}
\label{punga}
Consider the Ramanujan series
$$P(q):=E_2(q):=1-24\sum_{m \geq 1}\left(\sum_{d|m}d \right)q^m$$
and assume $N$ is prime. Consider the series
$$g(q):=-\frac{1}{24}(P(q)-NP(q^N))=\frac{N-1}{24}+f(q)\in {\mathbb Z}_{(p)}[[q]],$$
where
\begin{equation}
\label{acestf}
f(q)=\sum_{m \geq 1}\left(\sum_{A|m}\epsilon(A)A\right)q^m.\end{equation}
Then $g(q)$ is the Fourier expansion of a classical modular form in $M(\Gamma_0(N),{\mathbb Z}_{(p)},2)$
which is an eigenvector of the Hecke operators $T_2(n)$ for all $n \geq 1$ with eigenvalues $a_n:=\sum_{A|n}\epsilon(A)A$; cf. \cite{DI}, Example 2.2.6, Proposition 3.5.1, and Remark 3.5.2.
Let $\varphi:=\overline{f}=\sum_{m \geq 1} \overline{a}_m q^m \in q{\mathbb F}_p[[q]]$ be the reduction mod $p$ of $f(q)$.
By \cite{Knapp}, Theorem 9.17, the equalities (\ref{knacond}) hold with $\k=0$.
So by Theorem  \ref{nopiine}
the series
\begin{equation}
\label{inpat2}
f^{\sharp,2}=f^{\sharp,2}(q,q',q''):=\frac{1}{p}
\cdot \sum_{n \geq 1}\frac{a_n}{n}(\phi^2(q)^n-a_p \phi(q)^n+pq^n) \in {\mathbb Q}_p[[q,q',q'']]\end{equation}
belongs to ${\mathbb Z}_p[[q]][q',q'']\h$ and  its reduction mod $p$ equals
\begin{equation}
\label{baie2}
\overline{f^{\sharp,2}}:=\overline{f^{\sharp,2}(q,q',q'')}=\varphi^{(-1)}-\overline{a}_p V(\varphi)\frac{q'}{q^p}+   V^{2}(\varphi)\left(\frac{q'}{q^p}\right)^{p}
\in {\mathbb F}_p[[q]][q'].\end{equation}
Note also that $T_{2}(n)\varphi=\overline{a}_n \cdot \varphi$ for $(n,p)=1$ and $U\varphi=\overline{a}_p \cdot \varphi$. So by Theorem \ref{curcan}
   $\overline{f^{\sharp,2}}=\varphi^{\sharp,2}$ is
   an    eigenvector of the Hecke operators $nT_{0}(n)$, ``$pT_{0}(p)$", $(n,p)=1$, with eigenvalues $\overline{a}_n,\overline{a}_p$.
Also, by the same Theorem, if in addition $\overline{a}_p\neq 0$ and $\k=0$, then the series $\varphi^{\sharp,1}$ in (\ref{varphisharpunu})
   is also an eigenvector of the Hecke operators $nT_{\k}(n)$, ``$pT_{\k}(p)$", $(n,p)=1$, with eigenvalues $\overline{a}_n,\overline{a}_p$.
  Note that  if  $N \equiv 1$ mod $p$ then  Equations \ref{oprit} hold because
$$\begin{cases}
a_n=\sum_{A|n}\epsilon(A)A\equiv n \sum_{A|n}\epsilon(A)A^{-1},\ \ \text{mod $p$ for $(n,p)=1$},\\
\ \\
a_p=\sum_{A|p}\epsilon(A)A\equiv 1\ \ \text{mod $p$}.\end{cases}$$
Note also that if $N \equiv 1$ mod $p$ it follows that $f(q)\equiv g(q)$ mod $p$ so  $\varphi(q)$  is the Fourier expansion
  of a modular form in $M(\Gamma_0(N),{\mathbb F}_p,2)$
\end{example}

\subsection{Review of $\d$-modular forms \cite{eigen, igusa}}
Let $V$ be an affine smooth scheme over $R$ and fix a closed embedding $V \subset {\mathbb A}^m$ into an affine space over $R$.

\begin{definition}
A map $f:V(R)\ra R$ is called a {\it $\d$-function
of order $r$} on $X$ \cite{char} if
there exists a restricted power series $\Phi$ in $m(r+1)$ variables, with $R$-coefficients such that
$$f(a)=\Phi(a,\d a,...,\d^r a),$$
for all $a\in V(R)\subset R^m$. We denote by $\cO^r(V)$ the ring of $\d$-functions of order $r$ on $V$. \end{definition}

 (Recall that {\it restricted} means {\it with coefficients converging $p$-adically to $0$}; also the definition above does not depend on the embedding $V \subset {\mathbb A}^m$.)
 Composition with $\d$ defines $p$-derivations $\d:\cO^r(V)\ra \cO^{r+1}(V)$.
 The rings $\cO^r(V)$ have the following universality property: for any $R$-algebra homomorphism
 $u:\cO(V)\ra A$ where $A$ is a $p$-adically complete $\d$-ring there are unique $R$-algebra maps
 $u^r:\cO^r(V) \ra A$ that commute in the obvious sense with $\d$ and prolong $u$.

Let now $V$ be as in (\ref{ve}) with $B=R$ and ${\mathbb Z}[1/N,\zeta_N]\subset R$ a fixed embedding.

 \begin{definition} \cite{igusa}
A $\d-${\it modular function}  of order $r$ (on $\Gamma_1(N)$, holomorphic on $X$) is a $\d$-function
$f:V(R)\ra R$ of order $r$.\end{definition}

 Let $W:=\bZ[\phi]$ be the ring generated by $\phi$. For $w=\sum a_i\phi^i\in W$ ($a_i \in \bZ$)
 set $deg(w)=\sum a_i\in \bZ$; for $\lambda\in R^{\times}$ we set
 $\lambda^w:=\prod \phi^i(\lambda)^{a_i}$.

\begin{definition}
 A {\it $\d$-modular form of weight $w$} (of order $r$, on $\Gamma_1(N)$, holomorphic on $X$) is a $\d$-modular function $f:V(R)\ra R$ of order $r$ such that
 $$f(\lambda \cdot a)=\lambda^w f(a),$$
  for all $\lambda\in R^{\times}$ and $a\in V(R)$, where $(\lambda,a)\mapsto \lambda \cdot a$ is the natural action $R^{\times}\times V(R)\ra V(R)$.\end{definition}

 We  denote by  $M^r:=\cO^r(V)$ the ring of all $\d$-modular functions of order $r$ and we set $M^{\infty}:=\bigcup_{r \geq 0} M^r$.
We denote by $M^r(w)$ the $R$-module of $\d$-modular forms of order $r$ and weight $w$; cf. \cite{igusa}. (In \cite{eigen}
the space $M^r(w)$ was denoted by $M^r(\Gamma_1(N),R,w)$ or $M^r_{ord}(\Gamma_1(N),R,w)$ according as $X$ is $Y$ or $Y_{ord}$.) Note that $M^r(0)$ identifies with $\cO^r(X)$ which, in its turn, embeds into $\cO^r(X_1(N))$.

By the universality property of the rings $M^r=\cO^r(V)$ there exists a unique $\d$-ring homomorphism (the {\it $\d$-Fourier expansion map}) $$E:M^{\infty}\ra S^{\infty}_{for}:=\bigcup_{r\geq 0}R((q))[q',...,q^{(n)}]\h,\ \  \ E(f)=f(q,q',q'',...),$$
extending the Fourier expansion map $E:M\ra R((q))\h$.
  We may also consider the composition
 $$M^{\infty} \ra S^{\infty}_{for}\stackrel{\pi}{\ra} R((q))\h,\ \ \ f \mapsto f(q),$$
 where  the map $\pi$ sends $q',q'',...$ into $0$;
 we refer to this composition as the {\it Fourier expansion map}.

Recall the ``$\d$-expansion principle":

\begin{proposition}
\label{princip}
\cite{eigen}
The maps $E:M^r(w)\ra R((q))[q',...,q^{(r)}]\h$ are injective with torsion free cokernel;
hence the induced maps $\overline{E}:M^r(w)\otimes k \ra k((q))[q',...,q^{(r)}]$ are injective.
\end{proposition}

{\it Proof}.
This is \cite{eigen}, Lemma 6.1.
\qed

\medskip

Recall also the following  result:

\begin{theorem}
\label{zzz}
\cite{eigen}
If in Example \ref{gura} $\k=0$, $F={\mathbb Q}$, and $p\gg 0$  then the series $f^{\sharp,2}(q,q',q'')\in R[[q]][q',q'']\h$
in (\ref{inpat}) is the image
of a (unique)
 $\d$-modular form (still denoted by) $f^{\sharp,2}\in \cO^2(X_1(N)) \subset M^2(0)$.
 If in addition $f$ in Example \ref{gura} is of ``CL type" then the series $\varphi^{\sharp,1} \in k[[q]][q']$ in that Example  is the image of a $\d$-modular form $f^{\sharp,1}\in \cO^1(X_1(N))\subset M^1(0)$.
\end{theorem}

Here by $f$ being of {\it CL type} we mean that the Neron model of the elliptic curve over ${\mathbb Q}$ associated to $f$ via the Eichler-Shimura construction has good ordinary reduction and its base change to $R$ is the canonical lift of this reduction; cf. \cite{eigen,BP} for more details.

\medskip

{\it Proof}.
Let $f^{\sharp}\in \cO^r(X_1(N))$ be as in \cite{eigen},
Theorems 6.3 and 6.5; cf. also \cite{BP}, Lemma 4.18. So $r$ is $1$ or $2$ according as $f$ is or is not of CL type.
Then Theorem \ref{zzz} follows from  \cite{eigen},
Theorems 6.3 and 6.5, by letting the $\d$-modular form $f^{\sharp,2}$ be defined by
$$f^{\sharp,2}:=\begin{cases}
f^{\sharp},\ \ \text{if $f$ is not of CL type},\\
\phi(f^{\sharp})-a_pf^{\sharp},\ \ \ \text{if $f$ is of CL type},\end{cases}$$
and by letting
$$f^{\sharp,1}:=f^{\sharp}\ \ \ \text{if $f$ is of CL type}.$$
\qed

\begin{remark}
It is tempting to conjecture that if in Example \ref{gura} $\k\geq 0$ is arbitrary, $F={\mathbb Q}$,  and $p\gg 0$  then the series $f^{\sharp,2}(q,q',q'')$ is the $\d$-Fourier expansion
of a $\d$-modular form $f^{\sharp,2}\in M^r(\k)$ for some $r\geq 2$. An appropriate variant of this should also hold for arbitrary $F$. As we shall see, however, the situation is drastically different with Example
\ref{punga}; cf. Theorem \ref{contrast}.
\end{remark}

Recall the {\it Serre derivation operator}
$\partial: M\ra M$
introduced by Serre and Katz \cite{Katz}. (Cf. also \cite{book}, p.254 for a review). Recall that $\partial(L^{\otimes n}) \subset L^{\otimes (n+2)}$.  Recall also that if $X$ is contained in $Y_{ord}$ then one has the Ramanujan form $P \in M^0(2)$. By \cite{book}, Propositions 3.43, 3.45, 3.56, there exists a unique sequence of $R$-derivations $\partial_j:M^{\infty}\ra M^{\infty}$, $j \geq 0$, such that
\begin{equation}
\begin{cases}
\partial_j \circ \phi^s  =  0 \ \ \text{on $M$ for $j \neq s$}\\
\partial_j \circ \phi^j  =  p^j \cdot \phi^j \circ \partial\ \ \text{on $M$ for $j \geq 0$}\end{cases}
\end{equation}
These derivations then also have the property that
\begin{equation}
\begin{cases}
\partial_j = 0\ \ \text{on $M^{j-1}$ for $j \geq 1$}\\
\partial_j \circ \d^j=\phi^j \circ \partial\ \ \text{on $M$ for $j \geq 0$}\end{cases}
\end{equation}
and that
\begin{equation}
\label{pastreaza}
\partial_j(M^r(w)) \subset M^r(w+2\phi^j).
\end{equation}
Recall the Ramanujan theta operator $\theta=q\frac{d}{dq}:R((q))\ra R((q))$. Then by
\cite{book}, Lemma 4.18, there is a unique sequence of $R$-derivations $\theta_j:S^{\infty}_{for}\ra S^{\infty}_{for}$ such that
\begin{equation}
\label{conju}
\begin{cases}
\theta_j \circ \phi^s  =  0 \ \ \text{on $R((q))$ for $j \neq s$}\\
\theta_j \circ \phi^j  =  p^j \cdot \phi^j \circ \theta\ \ \text{on $R((q))$ for $j \geq 0$};\end{cases}
\end{equation}
and such that
\begin{equation}
\label{onem}
\begin{cases}
\theta_j = 0\ \ \text{on $R((q))[q',...,q^{(j-1)}]\h$ for $j \geq 1$}\\
\theta_j \circ \d^j=\phi^j \circ \theta\ \ \text{on $R((q))$ for $j \geq 0$.}\end{cases}
\end{equation}

\begin{proposition}
\label{barcfor}
For any $w=\sum_{i=0}^r a_i\phi^i\in W$, any $j \geq 0$, and any  $f \in M^r(w)$ the following formula holds in $S^{\infty}_{for}$:
$$E(\partial_j f)=\theta_j(E(f))-a_jp^jE(f)E(P)^{\phi^j}.$$
\end{proposition}

{\it Proof}.
This was proved in \cite{book}, Proposition 8.42 in the case of ``$\d$-Serre-Tate expansions"; the case of $\d$-Fourier expansions is entirely similar. (The level $1$ case of this Proposition was proved in \cite{Barcau} using the structure of the ring of modular forms of level $1$.)
\qed

Finally we recall the $\d$-modular forms $f^1$ and $f^{\partial}$ introduced in \cite{difmod} and \cite{Barcau} respectively:

\begin{proposition}
\label{funu}
\cite{difmod,Barcau,book} For each $r \geq 1$ there exists a unique form $f^r \in M^r(-1-\phi^r)$ such that
$$E(f^r)=\Psi^{\phi^{r-1}}+p\Psi^{\phi^{r-2}}+...+p^{r-1}\Psi.$$
In particular
$$E(f^r) \equiv \left(\frac{q'}{q^p}\right)^{p^{r-1}}\ \ \ mod\ \ \ p.$$
\end{proposition}

\begin{proposition}
\label{fpartz}
\cite{Barcau,book}
Assume $X=Y_{ord}$.
Then there exists a unique form $f^{\partial} \in M^1(\phi-1)$
such that
$E(f^{\partial})=1$. The form $f^{\partial}$ is invertible in the ring $M^1$ and its inverse belongs to $M^1(1-\phi)$.
Furthermore the image of $f^{\partial}$ in $M^1 \otimes k$, coincides with the image of the Eisenstein series $E_{p-1} \in M(\Gamma_1(N),R,p-1)$.
\end{proposition}

\begin{remark}
Note that Proposition \ref{funu} holds, in particular for $X=Y=Y_1(N)$. However Proposition \ref{fpartz}
fails for $X=Y$: the form $f^{\partial}$ has ``singularities" at the supersingular points.
\end{remark}

\subsection{Review of Katz' generalized $p$-adic modular forms \cite{Katz, Gouvea}}

Let $B$ be a $p$-adically complete ring, $p\geq 5$, and let $N$ be an integer coprime to $p$.
Consider the functor
\begin{equation}
\label{hite}
\{\text{$p$-adically complete $B$-algebras}\} \ra \{\text{sets}\}
\end{equation}
that attaches to any $A$ the set of isomorphism classes of triples $(E/A,\varphi, \alpha)$, where $E$ is an elliptic curve over $A$,
$\varphi$ is a trivialization, and $\alpha$ is an arithmetic level $\Gamma_1(N)$ structure. Recall that a {\it trivialization} is an isomorphism
between the formal group of $E$ and the formal group of the multiplicative group; an {\it arithmetic level $N$ structure}
is defined as an inclusion of flat group schemes over $B$, $\alpha: \mu_{N}\ra E$. So if $B$ contains a primitive $N$-th root of unity
(which we fix) then an arithmetic level $\Gamma_1(N)$ structure is the same as a level $\Gamma_1(N)$ structure. The functor (\ref{hite})
is representable by a $p$-adically complete ring $\bW(B,N)$.
The elements of this ring are called by Katz \cite{Katz} {\it generalized $p$-adic modular forms};  an element $f\in \bW(B,N)$ can be identified with a rule that naturally
attaches
to any test object $(E/A,\varphi,\alpha)$ an element $f(E,\varphi,\alpha)\in A$.
  Note that
 $\bW(B,N)=\bW(\bZ_p,N) \widehat{\otimes} B$.
Moreover there is a $\bZ_p^{\times}$-action on $\bW(B,N)$, $(\lambda,f) \mapsto \lambda \cdot f$,
 coming from the action of $\bZ_p^{\times}$ on the formal group of the multiplicative group.

There is a natural {\it Fourier expansion map}
$E:\bW(B,N) \ra B((q))\hat{\ }$ which is injective and has a flat cokernel over $B$ coming from evaluation on the Tate curve.
{\it From now on we shall view $\bW(B,N)$ as a subring of $B((q))\h$ via the Fourier expansion map.}

For $X=Y$ or $Y_{ord}$ note that the image of $\cO(V)=\bigoplus L^{\otimes \k} \ra R((q))\h$ is contained in $\bW$  and the morphism $\cO(V)\ra \bW$ is ${\mathbb Z}_p^{\times}$-equivariant with
$\lambda \in {\mathbb Z}_p^{\times}$ acting on $\eta \in L^{\otimes \k}$ via $(\lambda, \eta)\mapsto \lambda^{\k}\eta$.

Also $\bW(\bZ_p,N)$ possesses a natural ring endomorphism
$Frob$ which reduces modulo $p$ to the $p$-power Frobenius endomorphism of $\bW(\bZ_p,N) \otimes \bZ/p\bZ$. So if $R=\hat{\bZ}_p^{ur}$, as usual,
and if $\phi$ is the automorphism of $R$ lifting Frobenius then
$Frob \widehat{\otimes} \phi$ is a lift of Frobenius on
$$\bW:=\bW(R,N)=\bW(\bZ_p,N) \widehat{\otimes} R$$ which we denote  by $\phi_0$.
 Moreover the homomorphism $\bW(R,N) \ra R((q))\hat{\ }$ commutes with the action of $\phi_0$ where $\phi_0$  on $R((q))\hat{\ }$
is defined by $\phi_0(\sum a_nq^n):= \sum \phi(a_n)q^{np}$. Finally $\phi_0$ commutes with the action of ${\mathbb Z}_p^{\times}$.

Let $\chi:{\mathbb Z}_p^{\times} \ra {\mathbb Z}_p^{\times}$ be a continuous character.
An element $f \in \bW$ is said to have weight $\chi$ if $\lambda \cdot f=\chi(\lambda) f$ for all $\lambda\in {\mathbb Z}_p^{\times}$;
cf \cite{Serre, Gouvea}. We view integers $m \in {\mathbb Z}$ as identified with continuous characters by attaching to $m$ the character $\chi(\lambda)=\lambda^m$.
Recall from \cite{Gouvea}, p. 21 that the set of all $f \in \bW(B,N)\cap B[[q]]$ that have weight $\chi$ identifies with
 the set of {\it $p$-adic modular forms of weight $\chi$ defined over $B$} in the sense of Serre \cite{Serre}
i.e. the set of series in $B[[q]]$ which are $p$-adic limits of classical modular forms over $B$ of weights $\k_n\in {\mathbb Z}$ and level $N$ where $\k_n \ra \chi$.
Note that since $\phi_0$ commutes with the action of ${\mathbb Z}_p^{\times}$ on $\bW$ it follows that if $f \in \bW$ has weight $\chi$ then so does $\phi_0(f)\in \bW$.

\subsection{Application to $\d$-eigenforms}
\label{we}
As noted in \cite{igusa} the image of the Fourier expansion map $M^{\infty}\ra R((q))\h$ is contained in $\bW$; this is by the universality
property of $\cO^r(V)$ and by the fact that $\bW$ possesses a lift of Frobenius $\phi_0$ and hence it is naturally a $\d$-subring of $R((q))\h$.

\begin{proposition}
\label{nea}
The image of  $M^r(w)$ in $\bW$ consists of elements of  weight $deg(w)$.
\end{proposition}

{\it Proof}.
It is easy to see that one may replace $X$ in the statement above by an open set of it. So one may assume $L$ is free on $X$.
Let $x$ be a basis of $L$. Then any element $f \in M^r(w)$ can be written as $f=f_0\cdot x^w$ where $f_0 \in \cO^r(X)$. Now the image of $x$ in $\bW$ has weight $1$.
Since $\phi_0$
on $\bW$ preserves the elements of a given weight it follows that the image of $x^w$ in $\bW$ has weight $deg(w)$. On the other hand $f_0$ is a $p$-adic limit
of polynomials with $R$-coefficients in elements of the form $\d^i g_0$, where $g_0 \in \cO(X)$. Again, since $\phi_0$ sends elements of weight $0$
in $\bW$ into elements of weight $0$ the same is true for $\d:\bW \ra \bW$. Since the image of $g_0$ in $\bW$
has weight $0$ so does the image of $\d^i g_0$ in $\bW$ and hence so does the image of $f$.
\qed

Next we state our main applications to ``$\d$-eigenforms" (i.e. $\d$-modular forms whose $\d$-Fourier expansions are ``$\d$-eigenseries"). First we will prove:

\begin{theorem}
\label{dezastru}
Assume  $\overline{f}=\overline{f}(q,q')\in k[[q]][q']$  is not a $p$-th power in $k[[q]][q']$ and assume  $\overline{f}$ is
 the reduction mod $p$ of the $\d$-Fourier expansion of a $\d$-modular form in $M^r(w)$
  with $r\geq 0$, $\k: =deg(w)\geq 0$. Assume furthermore that
 $\overline{f}$ is an eigenvector of all Hecke operators
 $nT_{\k}(n)$, ``$pT_{\k}(p)$", $(n,p)=1$, with eigenvalues $\overline{\lambda}_n, \overline{\lambda}_p \in k$.
Then there exists  $\varphi=\varphi(q)\in qk[[q]]$
which is the Fourier expansion of a modular form
 in $M(\Gamma_1(N),k,\k')$, $\k'\geq 0$,
 $\k' \equiv \k+2$ mod $p-1$, and there exist
 $c, c_i \in k$, $i \geq 0$, with
 $\ee \cdot c_{i-1}=\overline{\lambda}_p c_i$ for $i \gg 0$,
  such that $\varphi$ is an eigenvector of all Hecke operators
 $T_{\k+2}(n)$, $T_{\k+2}(p)$, $(n,p)=1$, with the same eigenvalues $\overline{\lambda}_n, \overline{\lambda}_p$
  and such that
 $\overline{f}$ satisfies (\ref{zor}).
\end{theorem}

Conversely we will prove:

\begin{theorem}
\label{converseofdezastru}
Assume $\varphi \in qk[[q]]$ is the Fourier expansion of a  modular form
 in $M(\Gamma_1(N),k,\k')$, $\k' \geq 0$, $\k'\equiv \k+2$ mod $p-1$,
 and  that $\varphi$ is an eigenvector of all Hecke operators
  $T_{\k+2}(n)$, $T_{\k+2}(p)$, $(n,p)=1$, with eigenvalues
  $\overline{\lambda}_n, \overline{\lambda}_p \in k$. Assume $X=Y_{ord}$. Consider  the series
 $\overline{f}=\overline{f}(q,q')\in k[[q]][q']$ defined by the formula (\ref{zor}) with $c=0, c_i \in k$ for $i\geq 0$, and
 $c_i=0$ for $i \gg 0$.
 Then $\overline{f}$ is the $\d$-Fourier expansion of a $\d$-modular form
$f \in M^1(\k)$ and (by Theorem \ref{curcan}) is an eigenvector of all Hecke operators $nT_{\k}(n)$, ``$T_{\k}(p)$", $(n,p)=1$,
with the same eigenvalues $\overline{\lambda}_n, \overline{\lambda}_p$.
\end{theorem}

Note that Theorems \ref{dezastru} and \ref{converseofdezastru} imply Theorem \ref{maint} in the Introduction.
The one-to-one correspondence in Theorem \ref{maint} is given by $\varphi\mapsto \varphi^{\sharp,2}$ with $\varphi^{\sharp,2}$
defined by (\ref{zor}).

\medskip

{\it Proof of Theorem \ref{dezastru}}.
By Theorem \ref{arnab} all we have to show is that $\varphi$ in that Theorem is the Fourier expansion
of a modular form in $M(\Gamma_1(N),k,\k')$,
 $\k' \equiv \k+2$ mod $p-1$.
Since $\overline{f}$ is not a $p$-th power we may assume $c_0=1$.
Now if $\overline{f}(q,q')$ is the reduction mod $p$ of the $\d$-Fourier expansion
$$E(f)=f(q,q',...,q^{(r)})\in S^{\infty}_{for}$$
 of a $\d$-modular form $f \in M^r(w)$ then, by Proposition \ref{barcfor}, and Equations \ref{zor} and \ref{onem} we have
 the following congruences mod $p$ in $S^{\infty}_{for}$:
 $$\begin{array}{rcl}
 E(\partial_1 f) & \equiv & \theta_1(E(f))\\
 \  & \  & \  \\
\  & \equiv & -\overline{\lambda}_p V(\varphi)q^{-p}\theta_1(\d q)\\
\  & \  & \  \\
\  & \equiv & -\overline{\lambda}_p V(\varphi)q^{-p} \phi(\theta q)\\
\  & \  & \  \\
\  & \equiv & -\overline{\lambda}_p V(\varphi).
\end{array}$$

By Equation (\ref{pastreaza}) we have that $\partial_1 f\in M^r(w+2\phi)$.
So by Proposition \ref{nea} the image $E(\partial_1 f)(q,0,...,0)$ of $E(\partial_1 f)$ in $R((q))\h$
is an element of weight $\k+2$ in $\bW$. So
$E(\partial_1 f)(q,0,...,0)$
 is congruent mod $p$
to the Fourier expansion of a classical modular form of weight $\k' \equiv \k+2$ mod $p-1$.
So $\overline{\lambda}_p V(\varphi)$ is the Fourier expansion
of a modular form in $M(\Gamma_1(N),k,\k')$.

If $\overline{\lambda}_p \neq 0$ then $V(\varphi)$ is the Fourier expansion of a modular form in $M(\Gamma_1(N),k,\k')$ hence so is $\varphi=UV\varphi$ (because $U$ preserves the weight \cite{Gross}, p.458).

If $\overline{\lambda}_p =0$ then, by (\ref{defoffi}) we have $\varphi=\sum_{(n,p)=1}\overline{\lambda}_n q^n$ so $\varphi=\theta (\varphi^{(-1)})=\theta(\varphi_0)$. Now $\varphi_0$ is the image of $E(f)$ in $k[[q]]$ so, as above,   by Proposition \ref{nea}, $\varphi_0$ is the Fourier expansion of a modular form
in $M(\Gamma_1(N),k,\k'')$ where $\k'' \equiv \k$ mod $p-1$. But $\theta$ sends Fourier expansions of modular forms of weight $\k''$ into Fourier expansions of modular forms of weight $\k''+p+1$; cf. \cite{Gross}, p. 458. So $\varphi$ is the Fourier expansion of a modular form in
$M(\Gamma_1(N),k,\k''+p+1)$,  and we are done because $\k''+p+1 \equiv \k+2$ mod $p-1$.
\qed

\bigskip

{\it Proof of Theorem \ref{converseofdezastru}}.
Set $\k'=\k+2+(p-1)\nu$, $\nu\geq 0$.
Since $\varphi^{(-1)}(q)=\theta^{p-2} \varphi(q)$ by get that $\varphi^{(-1)}(q)$
 is the Fourier expansion of
a modular form  over $k$ of weight $\k'+(p-2)(p+1)=\k+(p-1)(p+\nu)$
hence $V^i(\varphi^{(-1)}(q))$ is the Fourier expansion of a modular form over $k$ of weight
$\k_{0,i}:=p^i(\k+(p-1)(p+\nu))$; the latter lifts to a modular form  $\Phi_{0,i}\in M(\Gamma_1(N),R,\k_{0,i})$ which can be viewed as an element in $M^0(\k_{0,i})$. Also $V^{i+1}(\varphi)$ and $V^{i+2}(\varphi)$ are Fourier expansions of modular forms over $k$ of weights $\k_{1,i}:=p^{i+1}\k'$ and $\k_{2,i}:=p^{i+2}\k'$ so they lift to modular forms  $\Phi_{i,1}\in  M(\Gamma_1(N),R,\k_{1,i})$ and $\Phi_{2,i}\in M(\Gamma_1(N),R,\k_{2,i})$ respectively. The latter can be viewed as elements  of $M^0(\k_{1,i})$ and $M^0(\k_{2,i})$ respectively.
Finally note that $f^1 \cdot f^{\partial}\in M^1(-2)$ and the Eisenstein form $E_{p-1}$ can be viewed as an element in $M^0(p-1)$; its inverse is an element in $M^0(1-p)$. Let $\lambda_p\in R$ be a lift of $\overline{\lambda}_p$. Note that $\k_{0,i}\equiv \k$ mod $p-1$; set $e_{0,i}:=\frac{\k-\k_{0,i}}{p-1}$.
Similarly $\k_{1,i}\equiv \k+2$ mod $p-1$ and $\k_{2,i}\equiv \k+2p$ mod $p-1$; set $e_{1,i}:=\frac{\k+2-\k_{1,i}}{p-1}$ and
$e_{2,i}:=\frac{\k+2p-\k_{2,i}}{p-1}$.
Then, by Propositions
\ref{funu} and \ref{fpartz}
$\overline{f}$ is the $\d$-Fourier expansion of the $\d$-modular form
\begin{equation}
\label{celalalt}
\sum_{i\geq 0} c_i\left[E_{p-1}^{e_{0,i}}\cdot \Phi_{0,i}-\lambda_p \cdot E_{p-1}^{e_{1,i}}\cdot \Phi_{1,i} \cdot  (f^1 \cdot f^{\partial})+\ee
\cdot E_{p-1}^{e_{2,i}} \cdot \Phi_{2,i} \cdot (f^1 \cdot  f^{\partial})^p\right]\end{equation}
which is an element of $M^1(\k).$ This ends the proof.
\qed

\begin{example}
\label{oua}
We consider a special case of Example \ref{gura}.
 Let
 \begin{equation}
 f(q)=\sum_{m\geq 1}a_mq^m\in q{\mathbb Z}[[q]]\end{equation} be the Fourier expansion of a
 cusp form  $f\in M(\Gamma_0(N),{\mathbb Z},2)$.
  Assume $a_1=1$ and assume $f(q)$  is an eigenvector for all the Hecke operators
 $T_{2}(n)$ with $n \geq 1$.
 Assume $p$ is a prime and let
  $\varphi:=\overline{f}=\sum_{m \geq 1} \overline{a}_m q^m \in q{\mathbb F}_p[[q]]$ be the reduction mod $p$ of $f(q)$.
Then the equalities (\ref{knacond}) hold with $\k=0$. So by Theorem  \ref{nopiine} the series
\begin{equation}
\label{inpatuletz}
f^{\sharp,2}=f^{\sharp,2}(q,q',q''):=\frac{1}{p}
\cdot \sum_{n \geq 1}\frac{a_n}{n}(p^{\k}\phi^2(q)^n-a_p \phi(q)^n+pq^n) \in {\mathbb Q}_p[[q,q',q'']]\end{equation}
belongs to ${\mathbb Z}_p[[q]][q',q'']\h$ and  its reduction mod $p$ equals
\begin{equation}
\label{baiemica}
\overline{f^{\sharp,2}}:=\overline{f^{\sharp,2}(q,q',q'')}=\varphi^{(-1)}-\overline{a}_p V(\varphi)\frac{q'}{q^p}+  V^{2}(\varphi)\left(\frac{q'}{q^p}\right)^{p}
\in {\mathbb F}_p[[q]][q'].\end{equation}
Note also that $T_{2}(n)\varphi=\overline{a}_n \cdot \varphi$ for $(n,p)=1$ and $U\varphi=\overline{a}_p \cdot \varphi$. So by Theorem \ref{curcan}
   $\overline{f^{\sharp,2}}$ is
   an eigenvector of the Hecke operators $nT_0(n)$, ``$pT_0(p)$", $(n,p)=1$, with eigenvalues
   $\overline{a}_n,\overline{a}_p$.
In addition, if $p\gg 0$, by Theorem \ref{zzz}, the series $f^{\sharp,2}(q,q',q'')$ in (\ref{inpatuletz}) is the $\d$-Fourier expansion of a $\d$-modular form $f^{\sharp,2}\in \cO^2(X_1(N))\subset M^2(0)$.

On the other hand, as in the proof, of Theorem \ref{converseofdezastru},
 $\varphi^{(-1)}(q)$ is the Fourier expansion of a modular form over $k$ of weight
$p^2-p$; the latter lifts to a modular form  $\Phi_{0}\in M(\Gamma_1(N),R,p^2-p)$ which can be viewed as an element in $M^0(p^2-p)$. Also $V(\varphi)$ and $V^{2}(\varphi)$ are Fourier expansions of modular forms over $k$ of weights $2p$ and $2p^2$ so they lift to modular forms  $\Phi_{1}\in  M(\Gamma_1(N),R,2p)$ and $\Phi_{2}\in M(\Gamma_1(N),R,2p^2)$ respectively. The latter can be viewed as elements  of $M^0(2p)$ and $M^0(2p^2)$ respectively.
Then
$\overline{f^{\sharp,2}(q,q',q'')}$ is the $\d$-Fourier expansion of the $\d$-modular form
\begin{equation}
\label{celalaltuletz}
f^{!}:=E_{p-1}^{-p}\cdot \Phi_{0}-a_p \cdot E_{p-1}^{-2}\cdot \Phi_1 \cdot  (f^1 \cdot f^{\partial})+
\cdot E_{p-1}^{-2p} \cdot \Phi_{2} \cdot (f^1 \cdot  f^{\partial})^p\in M^1(0).\end{equation}

Note now that $f^{\sharp,2}\in M^2(0)$ and $f^{!}\in M^1(0)$ have the same $\d$-Fourier expansion and the same weight.
By Proposition \ref{princip} (the ``$\d$-expansion principle") we get the following:
 \end{example}

 \begin{corollary}
 In the notation of Example \ref{oua}
 we have the congruence
$f^{\sharp,2}\equiv f^{!}$ mod $p$ in $M^2(0)$.
\end{corollary}

Note that the right hand side of this congruence has order $1$ and has a priori ``singularities" both at the cusps of $X_1(N)$ and at the supersingular points. In stark contrast with that, the left hand side of the above congruence has {\it no} ``singularity" at either the cusps or the supersingular points.

Also in  stark contrast with Theorem \ref{zzz}
 we have the following consequence of Theorem \ref{dezastru}.

\begin{theorem}
\label{contrast}
Let $f(q)$ be as in Example \ref{punga} and assume $N \not\equiv 1$ mod $p$ (for instance $p\gg 0$).
Then the series $\overline{f^{\sharp}(q,q',q'')}$ in (\ref{baie2}) is not the image of any element in any space $M^r(w)$ with $r \geq 0$, $deg(w)=0$.
\end{theorem}

{\it Proof}.
Assume  the notation of Example \ref{punga}. By Theorem \ref{dezastru} it follows that the image of $f(q)$
in ${\mathbb F}_p[[q]]$ is the Fourier expansion of some modular form $\widehat{f}\in M(\Gamma_1(N),{\mathbb F}_p,2+(p-1)\nu)$, $\nu\geq 0$. On the other hand, by Example \ref{punga} we know that the image of $g(q)$ in ${\mathbb F}_p[[q]]$ is the Fourier expansion of a modular form $\widehat{g}\in M(\Gamma_0(N),{\mathbb F}_p,2)$. It follows that the modular form
$$\widehat{h}:=E_{p-1}^{\nu}\cdot \widehat{g}-\widehat{f}\in M(\Gamma_1(N),{\mathbb F}_p,2+(p-1)\nu)$$
has Fourier expansion a constant  $\gamma:=\frac{N-1}{24} \in {\mathbb F}_p^{\times}$. On the other hand $\gamma$, viewed as an element in $M(\Gamma_0(N),{\mathbb F}_p,0)$ has Fourier expansion $\gamma$.
By the Serre and Swinnerton-Dyer Theorem \cite{Goren}, p.140, the difference $\widehat{h}-\gamma$
is divisible by $E_{p-1}-1$ in the ring $\bigoplus_{\k \in {\mathbb Z}}M(\Gamma_1(N),{\mathbb F}_p,\k)$.
It follows that the  weights
$2+(p-1)\nu$ and $0$ are congruent mod $p-1$,  a contradiction.
\qed

\bigskip

\bigskip
\bigskip

\bibliographystyle{amsplain}

\begin{thebibliography}{10}

\bibitem{Barcau} Barcau, M: Isogeny covariant differential
modular forms and the space of elliptic curves up to isogeny,
Compositio Math., 137 (2003), 237-273.

\bibitem{f1}
 Borger, J.,
$\Lambda$-rings and the field with one element, arXiv:math/0906.3146


\bibitem{char} Buium, A.:
Differential characters of Abelian varieties
     over $p-$adic fields, Invent. Math.
     122, 309-340 (1995).

\bibitem{pjets} Buium, A.: Geometry of $p$-jets, Duke Math. J., 82, 2, (1996), 349-367.

\bibitem{difmod} Buium, A.: Differential modular forms,
 Crelle J., 520 (2000), 95-167.

\bibitem{book} Buium, A.:
Arithmetic Differential Equations. Math. Surveys and Monographs
118, AMS (2005)

\bibitem{dcc} Buium, A.: Differential characters on curves,
Serge Lang memorial volume, to appear



\bibitem{eigen} Buium A., Differential eigenforms, J. Number Theory 128 (2008), 979-1010.

\bibitem{igusa} Buium A., Saha A., The ring of differential Fourier expansions, preprint.

\bibitem{BP} Buium, A., Poonen, B.:
Independence of points on elliptic curves arising from special points
on modular and Shimura curves, II: local results,
Compositio Math., 145 (2009), 566-602.

\bibitem{DR} Deligne, P., Rappoport, M., Schemas de modules de courbes elliptiques, LNM 349, Springer 1973, pp. 143-316.

\bibitem{DI} Diamond, F., and Im, J.:
Modular forms and modular curves. In:Seminar on Fermat's Last
Theorem, Conference Proceedings, Volume 17, Canadian Mathematical
Society, pp. 39-134 (1995).


\bibitem{Goren} Goren, E. Z., Lectures on Hilbert Modular Varieties and Modular Forms, CRM Monograph Series CRMM 14, 2002.

\bibitem{Gouvea} Gouvea, F., Arithmetic of $p$-adic modular forms, Lecture Notes in Math. 1304, Springer, 1985.

    \bibitem{Gross} Gross, B. H., A tameness criterion
for Galois representations associated to modular forms mod $p$,
Duke Math. J., 61, 2, 445-517 (1990)

\bibitem{Hida} Hida, H., Geometric modular forms and elliptic curves, World Scientific (2000).


\bibitem{Katz}  Katz, N.: $p-$adic properties of
modular schemes and modular forms, LNM 350, Springer, Heidelberg
(1973).

\bibitem{Knapp} Knapp, A.: Elliptic Curves, Math. Notes,
Princeton Univ. Press (1992)

\bibitem{Serre} Serre, J. P.: Formes modulaires et fonctions
z\'{e}ta $p-$adiques. In: LNM 350 (1973).



\end{thebibliography}

\end{document}